\documentclass[12pt]{article}

\usepackage{graphicx}
\usepackage{latexsym,amsmath,amsfonts,amscd, amsthm}
\usepackage{color}
\usepackage{bm}
\usepackage{tikz}
\usepackage{multirow}
\usetikzlibrary{arrows,backgrounds,snakes}
\usepackage{epsfig,subfigure}

\newtheoremstyle{plainNoItalics}{}{}{\normalfont}{}{\bfseries}{.}{ }{}

\theoremstyle{plain}
\newtheorem{thm}{Theorem}[section]

\theoremstyle{plainNoItalics}

\newtheorem{defn}[thm]{Definition}
\newtheorem{rem}[thm]{Remark}
\newtheorem{prop}[thm]{Proposition}

\newcommand{\f}{\frac}

\newcommand{\beq}{\begin{equation}}
\newcommand{\eeq}{\end{equation}}
\newcommand{\beqa}{\begin{eqnarray}}
\newcommand{\eeqa}{\end{eqnarray}}
\newcommand{\bit}{\begin{itemize}}
\newcommand{\eit}{\end{itemize}}
\newcommand{\bedef}{\begin{defn}}
\newcommand{\edefn}{\end{defn}}
\newcommand{\bpro}{\begin{prop}}
\newcommand{\epro}{\end{prop}}

\newcommand{\df}{\partial}

\newcommand{\Dt}{\Delta t}

\newcommand{\eps}{\varepsilon}
\newcommand{\mD}{{\mathcal D}}

\newcommand{\xL}{{x_{i-\frac{1}{2}}}}
\newcommand{\xR}{{x_{i+\frac{1}{2}}}}
\newcommand{\iL}{{i-\frac{1}{2}}}
\newcommand{\iR}{{i+\frac{1}{2}}}

\begin{document}




\begin{center}
{\bf
High Order {Hierarchical} Asymptotic Preserving Nodal Discontinuous Galerkin IMEX Schemes For The BGK Equation \footnote{Research supported by NSF DMS-1217008, DMS-1522777 and Air Force Office of Scientific Computing FA9550-12-0318 and the Fundamental Research Funds for the Central Universities No. 20720160009.}
}
\end{center}

\vspace{.2in}
\centerline{
Tao Xiong \footnote{{School of Mathematical Sciences, Fujian Provincial Key Laboratory of Mathematical Modeling and High-Performance Scientific Computing, Xiamen University, Xiamen, Fujian, P.R. China, 361005. Email: txiong@xmu.edu.cn}}
and
Jing-Mei Qiu\footnote{Department of Mathematics, University of Houston, Houston, 77004. E-mail: jingqiu@math.uh.edu. }
}

\bigskip
\noindent{\bf Abstract}
A class of high order asymptotic preserving (AP) schemes has been developed
for the BGK equation in Xiong et.~al.~(2015) \cite{JLQX_BGK}, which is based on the micro-macro formulation of the equation. The nodal discontinuous Galerkin (NDG) method with Lagrangian basis functions for spatial discretization and globally stiffly accurate implicit-explicit (IMEX) Runge-Kutta (RK) scheme as time discretization are introduced with asymptotic preserving properties. However, it is only necessary to solve the kinetic equation when the hydrodynamic description breaks down. Motivated by the recent work in Filbet and Rey (2015) \cite{Filbet_hybrid}, it is more naturally to construct a hierarchy scheme under the NDG-IMEX framework without hybridization, as the formal analysis in \cite{JLQX_BGK} shows that when $\eps$ is small, the NDG-IMEX scheme becomes a local discontinuous Galerkin (LDG) scheme for the compressible Navier-Stokes equations, and when $\eps=0$ it is a discontinuous Galerkin (DG) scheme for the compressible Euler equations. Moveover, we propose to combine the kinetic regime with the hydrodynamic regime including both the compressible Euler and Navier-Stokes equations. Numerical experiments demonstrate very decent performance of the new approach. In our numerics, all three regimes are clearly divided, leading to great savings in terms of the computational cost.

\bigskip
\noindent {\bf Keywords:}
Hierarchy scheme,
Compressible Euler Equations,
Compressible Navier-Stokes Equations,
Asymptotic Preserving,
Nodal Discontinuous Galerkin,
IMEX,
BGK equation

\section{Introduction}
\label{sec1}
\setcounter{equation}{0}
\setcounter{figure}{0}
\setcounter{table}{0}

In physics, rarefied gases can be modeled by kinetic description using the Boltzmann equation. In such a description, Knudsen number $\eps$ is an important dimensionless parameter, defined as the ratio of the molecular mean free path length to a representative physical length scale, characterizing the frequency of molecular collisions or how rarefied the gas is. In the zero limit of Knudsen number, the compressible Euler system describing the conservation of mass, moment and energy is a sufficient macroscopic model, while when the Knudsen is sufficiently small but not zero, the compressible Navier-Stokes equations including a correction term on viscosity and heat conductivity are needed. BGK equation is a simplified model for the Boltzmann equation, which is introduced by Bhatnagar, Gross and Krook \cite{bhatnagar1954model}, in a hyperbolic scaling.

Many numerical schemes have been proposed for solving the BGK and Boltzmann equations with a wide range of Knudsen number. A micro-macro decomposition framework was proposed by Bennoune, Lemou, Mieussen \cite{bennoune2008uniformly}, which can correctly capture the macroscopic Navier-Stokes limit when the Knudsen number is sufficiently small. Various versions of implicit-explicit schemes were proposed for the BGK equations in  \cite{pieraccini2007implicit, pieraccini2012microscopically} and for the ES-BGK equation in \cite{filbet2011asymptotic}. A BGK-penalization strategy was proposed by Filbet and Jin \cite{filbet2010class} for the Boltzmann equation. These methods are all related to the asymptotic preserving (AP) schemes, which are designed to mimic the asymptotic limit from the kinetic to the hydrodynamic models on the PDE level as $\eps$ goes to 0 \cite{jin2010asymptotic}.

A family of high order AP schemes for the BGK equation has been developed in \cite{JLQX_BGK}, based on the micro-macro decomposition framework. The proposed methods work for both constant Knudsen number $\eps$ and spatially variant $\eps=\eps(x)$ in a wide range. The high order spatial accuracy is achieved by nodal discontinuous Galerkin (NDG) finite element approaches \cite{hesthaven2008nodal}, and the high order temporal accuracy is achieved by globally stiffly accurate implicit-explicit (IMEX) Runge-Kutta (RK) methods \cite{boscarino2011implicit,boscarino2010class}. A formal asymptotic analysis showing that the scheme becomes a DG method \cite{shu2003high} with explicit RK time discretizations for the compressible Euler system in the zero limit of the Knudsen number. While for sufficiently small $\eps$ it gives rise to a local DG (LDG) discretization \cite{bassi1997high, cockburn1998local, baumann1999discontinuous, lomtev1999discontinuous, bassi2002numerical}, up to $\mathcal{O}(\eps^2)$, for the compressible Navier-Stokes equations.

Although it is more accurate to use kinetic models to describe physics problems, it is computationally very expensive to simulate. On the other hand, fluid descriptions, such as compressible Euler and Navier-Stokes equations, typically break down near shocks or kinetic boundary layers. In a multi-scale scenario, it is of interest to use the kinetic model only locally in regions where it is necessary, while taking the advantage of low computational cost of the fluid system elsewhere. For computational efficiency, many hybrid kinetic/fluid schemes with automatic domain decomposition criteria have been developed. Many of these criteria are based on the macroscopic quantities to pass from the hydrodynamic description to kinetic ones. They are easy to compute numerically, but they could become inaccurate near shock or boundary layers. For example, Boyd, Chen and Candler \cite{Boyd1995} proposed a criterion based on the local Knudsen number, where the kinetic description is used when the quantity is below a problem-dependent threshold value. This criterion later was used by Kolobov et al.
with a discrete-velocity model of the Boltzmann equation and a kinetic scheme for the hydrodynamic equations \cite{Kolobov2007}, and then by Degond and Dimarco with a Monte-Carlo solver for the kinetic equation and a finite volume method for the macroscopic ones.
Another criterion based on the viscous and heat fluxes of the Navier-Stokes equations, through
a Grad's 13-moments expansion was introduced by Tiwari in \cite{Tiwari1998}. This criterion
is used with a deterministic solver for the kinetic one by Degond, Dimarco and Mieussens in \cite{Degond2010}, Tiwari, Klar and Hardt in \cite{Tiwari2009, Tiwari2012}, Alaia and Puppo in \cite{Puppo2012} and Dimarco, Mieussens and Rispoli in \cite{Dimarco2014}. Recently in \cite{Filbet_hybrid} Filbet and Rey proposed a hybrid method based on the moment realizability criteria introduced by Levermore, Morokoff and Nadiga \cite{David_moment}. In this work, the criteria to/from kinetic from/to hydrodynamic regimes via macroscopic and microscopic quantities are proposed respectively. The hybrid scheme combines a central finite volume scheme using central Lax-Friedrichs fluxes \cite{Nessyahu1990} for the fluid equations with an asymptotic scheme with a first order IMEX discretization \cite{filbet2011asymptotic} for the kinetic ES-BGK equation.

Motivated by the criteria developed in \cite{Filbet_hybrid}, in this paper, we design a hierarchy scheme
based on the NDG-IMEX developed in \cite{JLQX_BGK}. The domain decomposition approach can be very
applied to the NDG-IMEX method naturally, as the scheme automatically becomes a fluid
solver in the hydrodynamic regime (a DG scheme for the compressible Euler equations and an LDG scheme for the compressible Navier-Stokes equations). Moreover, as a new ingredient, we propose a criterion to adaptively identify the Euler, Navier-Stokes and kinetic regimes, in which the corresponding high order numerical solvers are applied. Numerical experiments on one dimensional problems are performed to showcase the effectiveness of the new approach. Significant savings on the computational cost are observed, as compared to the full NDG-IMEX scheme for the kinetic BGK equation.

The rest of the paper is organized as follows. In Section 2, the BGK equation in a hyperbolic scaling and its micro-macro decomposition is given. In Section 3, high order AP nodal DG spatial discretization and globally stiffly accurate IMEX temporal discretizations are presented. The regime indicators are introduced. In Section 4, numerical results are performed for one dimensional problems. Conclusions are given in the final section.

\section{BGK Equation and Macro-micro Formulation}
\label{sec2}
\setcounter{equation}{0}
\setcounter{figure}{0}
\setcounter{table}{0}

\newcommand \smu {\sqrt{\mu}}
\newcommand \del {\partial}
\newcommand \ep {\varepsilon}


We consider the BGK equation in a hyperbolic scaling:
\begin{equation}
\partial_{t}f + v\cdot \nabla _{x}f =  \frac{1}{\ep}(M_U-f)
\label{bgk}
\end{equation}
where $f=f(x,v,t)$ is the distribution function of particles that depends on time $t>0$, position $x\in\Omega_x\subset \mathbb{R}^d$ and velocity $v\in \mathbb{R}^d$ for $d\ge 1$.
The parameter $\ep$ is the Knudsen number proportional to the mean free path, and $M_U$ is  the local Maxwellian defined by
\begin{equation}
M_U=M_U(x,v,t)=\frac{\rho(x,t)}{(2\pi T(x,t))^{d/2}}
\exp \left(-\frac{|v-u(x,t)|^2}{2T(x,t)} \right).
\label{maxwellian}
\end{equation}
$\rho$, $u$, $T$ represent the macroscopic density, the mean velocity, and the temperature respectively. $U$ has the components of the density, momentum and energy, which
are obtained by taking the first few moments of $f$:
\begin{equation}\label{U-vector}
U:=\left(
\rho,
\rho u,
E
\right)^\top = \int_{\mathbb{R}^d} \left(
1,
v,
\frac12 |v|^2
\right)^\top   f(v)dv.
\end{equation}
where $E=\frac12\rho|u|^2 +\frac{d}{2}\rho T$ and the superscript $\top$ denotes the transpose of the corresponding vector.
In this paper, we use $m=m(v):= \left( 1,v, \frac12 |v|^2 \right)^\top$ and let $\langle g \rangle: = \int_{\mathbb{R}^d} g(v) dv$. It is easy to check that $\langle m M_U  \rangle = \left(
\rho,
\rho u,
\frac12\rho|u|^2 +\frac{d}{2}\rho T
\right)^\top = U $. Hence $\langle  m (M_U-f) \rangle =0$, namely the BGK operator satisfies the conservation of mass, momentum and energy. Moreover, it enjoys the entropy dissipation: $\langle (M_U-f) \log f\rangle \leq 0$.

In the following, we briefly recall the micro-macro decomposition of \eqref{bgk}, from which the compressible Euler and Navier-Stokes limits will be followed. For details, see \cite{JLQX_BGK}.
Let us first introduce several notations.
Taking $M=M_U$ for short, we use $L^2_M$ to denote the Hilbert space equipped with the weighted inner product
\[
(f,g)_M:= \langle f  g M^{-1} \rangle,
\]
then for any function $f\in L^2_M$, we can write $f\in L^2_M$ as
\[
f= \Pi_M f +(\mathbf{I}-\Pi_M)f.
\]
where ${\Pi_M}f$ is an orthogonal projection from $L^2_M$ onto
$\mathcal{N}:=\text{span}\,\{M,vM,|v|^2M\}$, and it is explicitly given by
\begin{equation}\label{pi}
 {\Pi_M}f = \left( \frac1\rho\langle f\rangle +\frac{\langle (v-u)f\rangle}{\rho T}\cdot(v-u) +\frac{2}{d\rho} \langle (\frac{|v-u|^2}{2T}-\frac{d}{2}) f \rangle (\frac{|v-u|^2}{2T}-\frac{d}{2})  \right) M.
\end{equation}
By the orthogonal project, we can decompose $f$ into a macroscopic part $M$ and a microscopic part $\ep g$,
\begin{equation}\label{mm}
f := M+\ep g,
\end{equation}
with $\langle m g \rangle =0$. Inserting \eqref{mm} into \eqref{bgk} and applying the orthogonal projections ${\Pi_{M}}$ and $\mathbf{I}-\Pi_{M}$ respectively,
we will have the following macro-micro decomposed equations:
\begin{subequations}
\label{eq:mmdc}
\begin{align}
&\partial_t U + \nabla_x\!\cdot\!  F(U)+\eps  \nabla_x\!\cdot\!  \langle vmg \rangle=0, \label{equ}\\
&\eps \del_t g +  \eps (\mathbf{I}-\Pi_{M}) (v\!\cdot\!\nabla_x g)
=- \big( g +(\mathbf{I}-\Pi_{M})  (v\!\cdot\!\nabla_x M)\big). \label{eqg}
\end{align}
\end{subequations}
where the flux $F(U)=(\rho u, \rho u\otimes u +p I, (E+p)u)^\top$ and $p=\rho T$.
$I$ is the $d\times d$ identity matrix.
In a more general setting where the Knudsen number depends on the position $\ep=\ep(x)$,
the micro-macro formulation \eqref{eq:mmdc} should be written as follows:
\begin{subequations}
\label{eq:mmd}
\begin{align}
&\partial_t U +\nabla_x\!\cdot\!  F(U) + \nabla_x \!\cdot\! \big( \eps(x)  \langle vmg \rangle \big)=0, \label{equx} \\
&\eps(x)\partial_t g + ({\mathbf I}-\Pi_M)\nabla_x \!\cdot\!\big(\eps(x)vg\big)=-\big( g+({\mathbf  I}-\Pi_M)(v \!\cdot\!  \nabla_x M) \big). \label{eqgx}
\end{align}
\end{subequations}

We observe that the first two terms in \eqref{equ} form the Euler system and eq.~\eqref{equ} formally converges to the Euler system as $\ep\rightarrow 0$. For the third term in eq.~\eqref{equ}, the leading order ($\ep$ term) will give rise to the viscous term in compressible Navier-Stokes equations. To see this, we have, 
from \eqref{eqg}, 
\begin{equation}
g = -(\mathbf{I}-\Pi_{M}) (v\!\cdot\!\nabla_x M) +{\mathcal O}(\ep)
\end{equation}
and the direct computation shows that
\begin{equation}
\label{eq: analytic_deri}
\frac{(\mathbf{I}-\Pi_{M}) (v\!\cdot\!\nabla_x M) }{M}= \frac12 A: \left( \nabla_xu+(\nabla_xu)^\top -\frac{2}{d}(\nabla_x\!\cdot\! u)I \right) + B\cdot \frac{\nabla_xT}{\sqrt{T}}
\end{equation}
where
\beq
\label{eq: AB}
A=\frac{(v-u)\otimes(v-u)}{T} -\frac{|v-u|^2}{dT}I \;\text{ and }\;
B=\left(  \frac{|v-u|^2}{2T} - \frac{d+2}{2} \right)\frac{v-u}{\sqrt{T}}.
\eeq
Therefore, we deduce that
\begin{equation}\label{AB}
g =-  A: \left( \nabla_xu+(\nabla_xu)^\top -\frac{2}{d}(\nabla_x\!\cdot\! u)I \right) M - B\cdot \frac{\nabla_xT}{\sqrt{T}} M +{\mathcal O}(\ep).
\end{equation}
Here $X:Y = \sum_{i,j}X_{ij}Y_{ij}$ is the Frobenius inner product for matrices.
As we insert the expression \eqref{AB} into \eqref{equ}, we obtain
\begin{equation}
\begin{split}
\del_t \left(
\begin{array}{c} \rho\\ \rho u \\ E
\end{array}
 \right)+\nabla_x\!\cdot\! \left(
 \begin{array}{c}
 \rho u \\ \rho u\otimes u +p I  \\ (E+p)u
 \end{array}
 \right) = \ep  \left(
 \begin{array}{c} 0 \\ \nabla_x\!\cdot\!\sigma   \\ \nabla_x\!\cdot\!(\sigma u+ q)
\end{array}
 \right) +{\mathcal O}(\ep^2)
\end{split}\label{cns-bgk2}
\end{equation}
where
\beq
\label{eq: sigma_q}
\sigma = \mu D(u), \quad D(u)=\nabla_xu+(\nabla_xu)^\top -\frac{2}{d}(\nabla_x\!\cdot\! u)I \text{ and } q=\kappa \nabla_x T
\eeq
and
\[
\mu = \frac{T}{4} \langle A : A M \rangle \text{ and }\kappa = T\langle B \cdot B M \rangle.
\]
We refer to \cite{bardos1991fluid} for more details on the derivation. The above system \eqref{cns-bgk2} is the compressible Navier-Stokes equations if we disregard high order terms ${\mathcal O}(\ep^2)$. We note that when $d=1$, $\sigma=0$ and $\kappa=\frac{3}{2}\rho T$.

As we point out in \cite{JLQX_BGK}, although the BGK equation shares the basic properties of hydrodynamics with the Boltzmann equation, the Navier-Stokes equations derived from those equations display different viscosity and heat conductivity coefficients {\cite{cercignani1969mathematical, cercignani1988boltzmann, cercignani2000rarefied}.
We would remark that in this case, the ellipsoidal statistical BGK (ES-BGK) \cite{ES-BGK} operator
can be used {in the macro-micro decomposition framework}.


\section{NDG-IMEX Methods}
\label{sec:NDG-IMEX:3}

\setcounter{equation}{0}
\setcounter{figure}{0}
\setcounter{table}{0}

In this section, we will briefly review the NDG-IMEX scheme developed in \cite{JLQX_BGK}, where the nodal discontinuous Galerkin (NDG) method in space together with implicit-explicit (IMEX) Runge-Kutta (RK) time discretization is used. We will only focus on the one-dimensional case with $d=1$, $\Omega_x=[a, b]$ and $\Omega_v=[-V_c, V_c]$ with $V_c$ sufficiently large so that the Maxwellian defined in \eqref{maxwellian} can be regarded as zero outside $\Omega_v$ numerically. For simplicity, we will just consider Scheme II in \cite{JLQX_BGK} for the general $\ep(x)$ in \eqref{eq:mmd} here, but Scheme I can be used similarly. Extension to high dimensions will be considered later and contribute to our future work.

\subsection{Semi-discrete NDG methods}
\label{sec:3.1}

Start with a partition of $\Omega_x$, $a=x_{\frac12}<x_{\frac{3}{2}}<\cdots <x_{N_x+\frac12}=b$. Let $I_i=[\xL, \xR]$ denote an element with its length $h_i$, and let $h=\max_{i=0}^{N_x}h_i$. Given any non-negative integer $K$, we define a finite dimensional discrete space,
\begin{equation}
Z_h^K=\left\{z\in L^2(\Omega_x): z|_{I_i}\in P^K(I_i), \forall i\right\},
\label{eq:DiscreteSpace}
\end{equation}
and its vector version is denoted as ${\bf Z}_h^K$. The local space $P^K(I)$ consists of polynomials of degree at most $K$ on $I$.
Note that functions in $Z_h^K$ are piecewise defined. To distinguish the left and right limits of a function $z\in Z_h^K$ at a grid point $x_{i+\frac{1}{2}}$, we let
$z_{i+\frac{1}{2}}^\pm=\lim_{\Delta x\rightarrow \pm 0}z(x_{i+\frac{1}{2}}+\Delta x)$, and we also let $[z]_{i+\frac{1}{2}}=z_{i+\frac{1}{2}}^+-z_{i+\frac{1}{2}}^-$ as the jump.

Consider the relation \eqref{eq: analytic_deri} in one dimension, we have
\beq
\label{eq: analytic_deri:1d}
(\mathbf{I}-\Pi_{M})  (v\partial_x M)=A\frac{\partial_x T}{\sqrt{T}}M.
\eeq
With this, the equation \eqref{eqgx} becomes
\beq
\eps(x) \partial_t g + ({\mathbf I}-\Pi_M)\partial_x(\eps(x) vg)=-\left( g+A\f{\partial_x T}{\sqrt{T}} M \right). \label{eqgx2}
\eeq
A semi-discrete DG method for the micro-macro system \eqref{eq:mmd} is designed as following. We seek $U_h(\cdot,t)\in {\bf Z}_h^K$ and $g_h(\cdot,v,t)\in Z_h^K$, such that
$\forall \phi, \psi\in Z_h^K$ and $\forall i$,
\begin{subequations}
\label{eq:SDG}
\begin{align}
\int_{I_i} \partial_t U_h \phi dx=&\int_{I_i}\left( F(U_h)+\eps(x) \langle vm g_h\rangle\right)  \frac{d \phi(x)}{dx} dx - \hat{F}_\iR \phi^-_\iR
+ \hat{F}_\iL \phi^+_\iL   \notag \\
&-\eps(x_\iR)  \widehat{\langle vmg_h\rangle}_\iR \phi_\iR^-+\eps(x_\iL) \widehat{\langle vmg_h\rangle}_\iL \phi^+_\iL,
\label{eq:SDG:a} \\
\int_{I_i} \eps(x) \partial_t g_h\psi dx+&
\int_{I_i} ({\mathbf I}-\Pi_{M_h})\left(\mD_{h,1}(\eps(x) vg_h)\right)\psi dx =- \int_{I_i}  g_h\psi dx
- \int_{I_i} A_h \frac{r_h}{\sqrt{T_h}} M_h \psi dx,\label{eq:SDG:c}
\end{align}
\end{subequations}
Here $M_h=M_{U_h}$ according to \eqref{maxwellian}. $\mD_{h,1}(\eps(x)vg_h)(\cdot,v,t)\in Z_h^K$ and $r_h\in Z_h^K$ are approximations of the spatial derivative of $\eps(x)vg$ and $T$, respectively, based on DG discretizations. Particularly, $\forall \psi\in Z_h^K$ and $\forall i$,
\begin{equation}
\int_{I_i} \mD_{h,1}(\eps(x) vg_h)\psi dx:=
-\int_{I_i} \eps(x) vg_h \frac{d \psi}{dx} dx+\eps(x_\iR) \widetilde{(vg_h)}_{\iR} \psi^-_\iR
-\eps(x_\iL) \widetilde{(vg_h)}_{\iL} \psi^+_\iL,
\label{eq:Dd1}
\end{equation}
where $\widetilde{vg}$ is an upwind numerical flux consistent to $vg$,
\beq
\label{eq:vg:upwind:L-1}
\widetilde{vg}:=
\left\{
\begin{array}{ll}
v g^-,&\mbox{if}\; v>0,\\
v g^+,&\mbox{if}\; v<0,
\end{array}
\right.
\eeq
and $\forall \varphi \in Z_h^K$ and $\forall i$
\beq
\label{eq: D3}
\int_{I_i} r_{h}\varphi dx = -\int_{I_i} T_h \frac{d \varphi}{dx} dx+\widehat{T}_{h,\iR} \varphi^-_\iR- \widehat{T}_{h,\iL} \varphi^+_\iL.
\eeq
Here $T_h$, a macroscopic quantity, and $A_h$ can be obtained from $U_h$ based on the energy $E$ under \eqref{U-vector} and \eqref{eq: AB} respectively.
The numerical flux $\hat{F}=\hat{F}(U_h^-, U_h^+)$ in \eqref{eq:SDG:a} is taken to be the global Lax-Friedrichs flux \cite{cockburn1989tvb2}. Here the subscripts $i\pm\frac{1}{2}$ are temporarily omitted for simplicity. We take the hatted fluxed $\widehat{ \langle vmg\rangle}$ and $\widehat{T}$ as the central fluxes $\widehat{\langle vmg\rangle} = {\langle vm(g^- + g^+)\rangle/2}, \widehat{T} = (T^+ + T^-)/2$, while the alternating right-left and central fluxes introduced in \cite{JLQX_BGK} can also be used.

We further use the nodal basis to represent functions in the discrete space $Z_h^K$, and approximating the integrals in the schemes by numerical quadratures.
Note that the discrete space $Z_h^K|_{I_i}$ is simply $P^K(I_i)$. We particularly choose the local nodal basis (also called Lagrangian basis) $\{\phi_i^k(x)\}_{k=0}^K$ associated with the $K+1$ Gaussian quadrature points $\{x^k_i\}^K_{k=0}$ on $I_i$, defined as below
\begin{equation}
\phi_i^k(x) \in P^K(I_i), \quad\textrm{and} \quad \phi_i^k (x_i^{k'})=\delta_{k k'},\quad k, k'=0,1,\cdots, K.
\end{equation}
Here $\delta_{k k'}$ is the Kronecker delta function. We further let $\{\omega_k\}^K_{k=0}$ denote the corresponding quadrature weights on the reference element
$(-\frac{1}{2}, \frac{1}{2})$.

With the nodal basis functions, \eqref{eq:SDG}-\eqref{eq:Dd1} and \eqref{eq: D3} can be equivalently stated with the test functions $\phi$, $\psi$ both being taken as $\phi_i^k, k=0, 1, \cdots, K$. We also replace all the integral terms in {\eqref{eq:SDG}-\eqref{eq:Dd1} and \eqref{eq: D3}} by their numerical integrations based on $(K+1)$-point Gaussian quadrature. The scheme now becomes: find
$U_h(\cdot,t)\in {\bf Z}_h^K$ and $g_h(\cdot,v,t)\in Z_h^K$, with
$U_h(x,t)|_{I_i}=\sum_{k=0}^K U_i^k(t) \phi_i^k(x)$, $g_h(x,v,t)|_{I_i}=\sum_{k=0}^K g_i^k(v,t)\phi_i^k(x)$, such that $\forall i, k$,
\begin{subequations}
\label{eq:SDG1}
\begin{align}
\omega_k h_i \frac{d U^k_i}{dt}=& \sum_{k'=0}^K \omega_{k'} h_i F(U^{k'}_i)\frac{d \phi^k_i(x)}{dx}\big |_{x=x^{k'}_i}-\hat F_\iR \phi_i^k(x_\iR^-)+\hat F_\iL \phi_i^k(x_\iL^+)\nonumber \\
&+\sum_{k'=0}^K \omega_{k'}h_i\eps(x^{k'}_i)\langle vm g^{k'}_i\rangle \frac{d \phi^k_i(x)}{dx}\big |_{x=x^{k'}_i}
-\eps(x_{\iR})\widehat{\langle vmg_h\rangle}_\iR \phi_i^k(x_\iR^-) \nonumber\\
&+\eps(x_\iL)\widehat{\langle vmg_h\rangle}_\iL \phi_i^k(x_\iL^+), \label{eq:SDG1:a} \\
\eps(x^k_i)\omega_k h_i \df_t g^k_i=&({\mathbf I}-\Pi_i^k)\Bigg(v\sum_{k'=0}^K \omega_{k'} h_i \eps(x^{k'}_i)g^{k'}_i \frac{d \phi^k_i(x)}{dx}\big |_{x=x^{k'}_i}
-\eps(x_{\iR})\widetilde {(vg_h)}_\iR \phi_i^k(x_\iR^-) \nonumber \\
\label{eq:SDG2:c}
&+\eps(x_{\iL})\widetilde{(vg_h)}_\iL \phi_i^k(x_\iL^+)\Bigg)- \omega_k h_i g^k_i  + A_i^k \omega_k h_i {r_i^k}M_i^k/{\sqrt{T_i^k}}, \\
\label{eq:SDG2:d}
\omega_k h_i r_i^k & = - \sum_{k'=0}^K \omega_{k'} h_i T_i^{k'} \frac{d \phi^k_i(x)}{dx}\big |_{x=x^{k'}_i}
+ \widehat{T}_{h,\iR} \phi_i^k(x_\iR^-) -\widehat{T}_{h,\iL} \phi_i^k(x_\iL^+).
\end{align}
\end{subequations}
Here $M_i^{k'}={M_h}|_{x=x^{k'}_i}$, $\Pi_i^k=\Pi_{M_i^k}$ and $r_i^k=r_h|_{x=x^{k}_i}$. And the nodal values of $T_i^k=T_h|_{x=x^{k}_i}$ and $A_i^k=A_h|_{x=x^{k}_i}$ are obtained from $U_i^k$ based on \eqref{U-vector} and \eqref{eq: AB}.

To the end, we also need to discretize the $v$-direction. In this work, $\Omega_v=[-V_c, V_c]$ is discretized uniformly with $N_v$ points,  $\{v_j\}_{j=1}^{N_v}$. For the integration in $v$,
the mid-point rule is applied, which is spectrally accurate for smooth functions with periodic
boundary conditions  or with a compact support \cite{boyd2001chebyshev}.
Such approach does not preserve the conservation properties of mass, moment and energy at the discrete level as in \cite{mieussens2000discrete}, yet in \cite{JLQX_BGK} we have found it is a sufficiently accurate discretization for all test cases that we have performed.

\subsection{IMEX time discretization}
\label{sec:3.2}

Now we will formulate the IMEX RK time discretizations for the semi-discrete schemes introduced in Section \ref{sec:3.1}. First we rewrite the scheme in a compact form as follows. Find $U_h(\cdot,t)\in {\bf Z}_h^K$, $g_h(\cdot,v,t), r_h(\cdot, t)\in Z_h^K$, such that
$\forall \phi, \psi, \varphi\in Z_h^K$ and $\forall i$,
\begin{subequations}
\label{eq: compact}
\begin{align}
(\partial_t U_h, \phi) + F_h(U_h, \phi) &= D_h( \eps(x) g_h, \phi), \label{eq: compact1}\\
(\eps(x) \partial_t g_h, \psi) + b_{h, v}(\eps(x)g_h, \psi) &= s^{(1)}_h(g_h, \psi) + s^{(2)}_{h, v}(U_h, r_h, \psi), \label{eq: compact2}\\
(r_h, \varphi) =  H_h(U_h, \varphi), \label{eq: compact3}
\end{align}
\end{subequations}
where
\begin{subequations}
\label{eq: S2}
\begin{align}
F_h(U_h, \phi) &= - \int_{\Omega_x} F(U_h)\frac{d \phi(x)}{dx} dx - \sum_i \hat{F}_{h, i+\frac12} [\phi]_{i+\frac12},\\
D_h( \eps(x) g_h, \phi) & =  \int_{\Omega_x} \eps(x) \langle vm g_h\rangle\frac{d \phi(x)}{dx} dx + \sum_i \eps(x_{i+\frac12}) \widehat{\langle vmg_h\rangle}_\iR [\phi]_{i+\frac12}, \\
b_{h, v}(\eps(x)g_h, \psi) &= \int_{\Omega_x} ({\mathbf I}-\Pi_{M_h}) \mD_{h,1}(\eps(x) vg_h)\psi dx,\\
s^{(1)}_h(g_h, \psi) &=-\int_{\Omega_x} g_h\psi dx, \quad
s^{(2)}_{h, v}(U_h, r_h, \psi) =  -\int_{\Omega_x} A_h \frac{r_h}{\sqrt{T_h}} M_h \psi dx, \label{eq: s1_q}\\
H_h(U_h, \varphi)&= {-\left(\int_{\Omega_x} T_h \frac{d\varphi}{dx} dx + \sum_i \hat{T}_{h, i+\frac12} [\varphi]_{i+\frac12}\right)\big |_{T_h=T_h(U_h)}.}
\end{align}
\end{subequations}

High order globally stiffly accurate IMEX schemes can be characterized by a double Butcher Tableau
\beq
\label{eq: B_table}
\begin{array}{c|c}
\tilde{c} & \tilde{A}\\
\hline
 & \tilde{b}^\top\end{array} \ \ \ \ \
\begin{array}{c|c}
{c} & {A}\\
\hline
 & {b^\top} \end{array},
\eeq
where $\tilde{A} = (\tilde{a}_{ij})$ is an $s\times s$  lower triangular matrix with zero diagonal for an explicit scheme, and $A = (a_{ij})$ is an $s\times s$ lower triangular matrix with the diagonal entries not all being zero for a diagonally implicit RK (DIRK) method. The coefficients $\tilde{c}$ and $c$ are given by the standard relations
\begin{eqnarray}\label{eq:candc}
\tilde{c}_i = \sum_{j=1}^{i-1} \tilde a_{ij}, \ \ \ c_i = \sum_{j=1}^{i} a_{ij},
\end{eqnarray}
and vectors $\tilde{b} = (\tilde{b}_j)$ and $b = (b_j)$ represent the quadrature weights for internal stages of the RK method. The IMEX RK scheme is defined to be {\em globally stiffly accurate} if $\tilde{c}_s = c_s = 1$ and {$a_{sj} = b_j$, $\tilde{a}_{sj} = \tilde{b}_j$}, $\forall j=1, \cdots, s$.

Now the fully-discrete scheme using the Butcher notation can be written as
follows. Given $U_h^n\in {\bf Z}_h^K$ and $g_h^n \in Z_h^K$, we look for $U_h^{n+1}\in {\bf Z}_h^K$ and $g_h^{n+1} \in Z_h^K$, such that $\forall \phi, \psi \in Z_h^K$,
\begin{subequations}
\label{eq: compact_d_rk}
\begin{align}
\left(U^{n+1}_h, \phi \right) &= \left(U^{n}_h, \phi \right) -\Dt \sum_{l=1}^{s} \tilde{b}_l \left(F_h(U^{(l)}_h, \phi) -D_h( \eps(x) g^{(l)}_h, \phi)\right), \label{eq: compact_d_rk1}\\
\left(\eps(x) g^{n+1}_h, \psi\right) & =\left(\eps(x) g^{n}_h, \psi\right)- \Dt \sum_{l=1}^{s} \tilde{b}_l b_{h, v}(\eps(x)g^{(l)}_h, \psi) + \Dt \sum_{l=1}^{s} b_l \left(s^{(1)}_h(g^{(l)}_h, \psi) + s^{(2)}_{h, v}(U^{(l)}_h, r^{(l)}_h, \psi)\right). \label{eq: compact_d_rk2}
\end{align}
\end{subequations}
Here the approximations at the internal stages of one RK step, $U_h^{(l)}\in {\bf Z}_h^K$ and $g_h^{(l)}, {r_h^{(l)}} \in Z_h^K$ with $l=1, \cdots, s$, satisfy
\begin{subequations}
\label{eq: compact_rk}
\begin{align}
\left(U^{(l)}_h, \phi \right) &= \left(U^{n}_h, \phi \right) -\Dt \sum_{j=1}^{l-1} \tilde{a}_{lj} \left(F_h(U^{(j)}_h, \phi) -D_h( \eps(x) g^{(j)}_h, \phi)\right), \label{eq: compact_rk1}\\
\left(\eps(x) g^{(l)}_h, \psi\right) & =\left(\eps(x) g^{n}_h, \psi\right)- \Dt \sum_{j=1}^{l-1} \tilde{a}_{lj} b_{h, v}(\eps(x)g^{(j)}_h, \psi) + \Dt \sum_{j=1}^{l} a_{lj} \left(s^{(1)}_h(g^{(j)}_h, \psi) + s^{(2)}_{h, v}(U^{(j)}_h, r^{(j)}_h, \psi)\right), \label{eq: compact_rk2}\\
{(r_h^{(l)}, \varphi)} &{= H_h(U_h^{(l)}, \varphi),} \label{eq: compact_rk3}
\end{align}
\end{subequations}
for any $\phi, \psi, \varphi \in Z_h^K$. One can solve the IMEX scheme in a stage-by-stage fashion for {$l=1, \cdots, s$}, that is, we first solve $U^{(l)}_h$ explicitly from the equation \eqref{eq: compact_rk1}, then plug $U^{(l)}_h$ into \eqref{eq: compact_rk3} to solve $r_h^{(l)}$, and finally solve $g^{(l)}_h$ from \eqref{eq: compact_rk2}.

The third order IMEX scheme we use in our simulations is the globally stiffly accurate ARS(4, 4, 3) scheme \cite{ascher1997implicit} with a double Butcher Tableau
\beq
\label{eq: ars443}
\begin{array}{c|c c c c c}
0 & 0&0&0& 0 & 0\\
1/2 &1/2&0&0& 0 &0\\
2/3 &11/18&1/18&0& 0 &0\\
1/2 &5/6&-5/6&1/2& 0 &0\\
1 &1/4&7/4&3/4& -7/4 &0\\
\hline
&1/4&7/4&3/4& -7/4 &0\\
 \end{array} \ \ \ \ \
\begin{array}{c|c c c c c}
0 & 0&0&0& 0 & 0\\
1/2 &0&1/2&0&0 &0\\
2/3 &0&1/6&1/2& 0 &0\\
1/2 &0&-1/2&1/2& 1/2 &0\\
1 &0&3/2&-3/2& 1/2 &1/2\\
\hline
 &0&3/2&-3/2& 1/2 &1/2\\
 \end{array}
\eeq

For the fully discrete DG-IMEX scheme \eqref{eq: compact_d_rk}-\eqref{eq: compact_rk}, with operators specified in \eqref{eq: S2}, the Proposition 3.4 in \cite{JLQX_BGK} has shown that, for $0<\eps\ll 1$, the scheme is asymptotically equivalent, up to $\mathcal{O}(\eps^2)$,
to a local DG (LDG) method in its nodal form for the compressible Navier-Stokes equations
\begin{equation}
\begin{split}
\del_t \left(
\begin{array}{c} \rho\\ \rho u \\ E
\end{array}
 \right)+
 \del_x \left(
 \begin{array}{c}
 \rho u \\ \rho u^2 +p I  \\ (E+p)u
 \end{array}
 \right) = \ep  \del_x \left(
 \begin{array}{c} 0 \\
 0   \\
 \frac32 \rho T \partial_x T
\end{array}
 \right).
\end{split}\label{cns-bgk-1d}
\end{equation}
and the LDG scheme evolved in time by an explicit RK method characterized by a Butcher table $\tilde{A}$, $\tilde{b}$ and $\tilde{c}$ in \eqref{eq: B_table} is defined as follows: find $U^{n+1}_h(\cdot, t), U^{(l)}_h \in {\bf Z}_h^K$ and  $r^{(l)}_h(\cdot, t) \in Z_h^K$ with $l=1, \cdots, s$, such that $\forall \phi, \varphi \in Z_h^K$,
\beq
\label{eq: compact_cns_rk0}
\left(U^{n+1}_h, \phi \right) = \left(U^{n}_h, \phi \right) -\Dt \sum_{l=1}^{s} \tilde{b}_l \left(F_h(U^{(l)}_h, \phi) -\eps F^{(vis)}_h(U^{(l)}_h, r^{(l)}_h, \phi)\right),
\eeq
with
\beq
\label{eq: compact_cns_rk}
\left(U^{(l)}_h, \phi \right) = \left(U^{n}_h, \phi \right) -\Dt \sum_{j=1}^{l-1} \tilde{a}_{lj} \left(F_h(U^{(j)}_h, \phi) -\eps F^{(vis)}_h (U^{(j)}_h,  r^{(j)}_h, \phi)\right),
\eeq
and
\beq
\label{eq: q_rk}
(r_h^{(l)}, \varphi)=
-\sum_i\left(\int_{I_i} T_h^{(l)} \frac{d\varphi}{dx} dx + \widehat{T}^{(l)}_{h,i+\frac12} [\varphi]_{\iR}\right).
\eeq
Its nodal form can be similarly defined as in Section \ref{sec:3.1} and is omitted for brevity. When $\eps=0$, if we omit the $\eps$ terms in \eqref{eq: compact_cns_rk0} and
\eqref{eq: compact_cns_rk}, it becomes a RK DG scheme for the compressible Euler system.

\subsection{Regime indicators}
\label{sec:3.3}

In this section, we will introduce the regime indicators to group all the computational
cells into three classes: (I) Euler regime, (II) Navier-Stokes (NS) regime, (III) kinetic regime. We will start all cells in the kinetic regime unless the initial conditions are apparently in the fluid regime. We have the criteria for both directions: from/to kinetic to/from hydrodynamic regimes. In particular, we use macroscopic quantities to determine when the hydrodynamic description breaks down, and use microscopic ones to determine when the kinetic description is not necessary and a hydrodynamic description would be sufficient. Moreover, we also distinguish between compressible Euler and Navier-Stokes in the hydrodynamics regime,
where a finer criteria is proposed.

Before we start, we first follow \cite{Filbet_hybrid} to introduce several notations. For simplicity, in the following, we consider the problem in the one-dimensional case where $d=1$. We use the short hand notation for the rescaled microscopic velocity $V(v)=\frac{v-u}{\sqrt{T}},$
then $A$ and $B$ defined in \eqref{eq: AB} can be written as $A=0$ and $B(V)=\frac12\left(V^2-3\right)V$. We let
$\bar B :=\frac{1}{\rho}\int_{\mathbb{R}}B(V)f(v)dv.$
 


\subsubsection{From fluid to kinetic: The moment realizability criterion}
Following \cite{David_moment, Filbet_hybrid}, we define the moment realizability matrix as
\begin{equation}
\label{matrix}
{\bf M}:=\frac{1}{\rho}\int_{\mathbb{R}}{\bf m}\otimes {\bf m} f(v)dv,
\end{equation}
where ${\bf m}$ is the collisional invariant vector for $V=(v-u)/\sqrt{T}$,
\begin{equation}
{\bf m}:= \left(1, V, \frac{1}{\sqrt{2}}\left(V^2-1\right)\right).
\end{equation}
From the properties of the moments defined in \eqref{U-vector}, we have
\begin{align}
\label{mom}
{\bf M}&=\frac{1}{\rho}\int_{\mathbb{R}}
\left(\begin{array}{ccc}
1 & V & \frac{1}{\sqrt{2}}\left(V^2-1\right) \\
V & V^2 & \frac{1}{\sqrt{2}}\left(V^2-1\right)V \\
\frac{1}{\sqrt{2}}\left(V^2-1\right) & \frac{1}{\sqrt{2}}\left(V^2-1\right)V & \frac{1}{2}\left(V^2-1\right)^2
\end{array}\right)f(v)dv \nonumber \\
&=\left(\begin{array}{ccc}
1 & 0 & 0 \\
0 & 1 & \frac{1}{\sqrt{2}} \bar B \\
0 & \frac{1}{\sqrt{2}} \bar B & \bar C
\end{array}\right)
\end{align}
where $\bar C$ is the dimensionless fourth order moment of $f$:
\[
\bar C = \frac{1}{2\rho}\int_{\mathbb{R}}\left(V^2-1\right)^2 f(v)dv.
\]

Now let us consider different orders for the approximation
of $f$ with respect to $\eps$. The values of $\bar B$ and $\bar C$ can then be
explicitly determined.
\begin{itemize}
\item Zeroth order: Compressible Euler system. \\
We have $f=M_U$, since $\bar B$ only involves odd, centered moments of $f$, so
that $\bar B=0$. Direct computation shows that $\bar C=1$. In this case, ${\bf M}$
is the identity matrix, all the eigenvalues are $1$, we denote it as $\nu_{Euler}=1$.

\item First order: Compressible Navier-Stokes system. \\
We have $f=M_U+\eps g = M_U - \eps A \frac{T_x}{\sqrt{T}} M$,
with analytical expression for the projection term defined in \eqref{eq: analytic_deri}.
And we can find that $\bar B = -\eps \frac{\kappa}{\rho T^{3/2}} \nabla_x T$, where $q=\kappa \nabla_x T$ is
the heat flux. For $\bar C$, using symmetry arguments, we still have $\bar C=1$.
For the matrix ${\bf M}$, beside one eigenvalue to be $1$, the other two eigenvalues
would be $1\pm \bar B$. We denote the largest absolute eigenvalue for the
compressible Navier-Stokes equations as $\nu_{NS}=1+|\bar B|=1+\eps \frac{\kappa}{\rho T^{3/2}} |\nabla_x T|$.

\item Second order: Burnett equations. \\
The second order Burnett equations will be used as a reference for the kinetic equation, when the fluid is far away from the thermal equilibrium. It would be complicate to derive
the explicit expression for $\bar B$ under the micro-macro decomposition framework up to second order. Instead, we directly use the one obtained from Chapman-Enskog expansion as in \cite{Filbet_hybrid}. Since $D(u)=0$, here the expression for $\bar B$ would be
\beq
\bar B = -\eps \frac{\kappa}{\rho T^{3/2}} T_x - \eps^2 \frac{\mu^2}{\sqrt{T}} \left(
\frac{25}{6}u_x T-\frac{5}{3} \left(T u_{xx} + 7u_x T_x\right)\right).
\label{burnettB}
\eeq
Similarly, beside $1$ is one eigenvalue of ${\bf M}$, the other two eigenvalues are
$1\pm \bar B$, so that we denote $\nu_B=1+|\bar B|$ with $\bar B$ in \eqref{burnettB}.
\end{itemize}

A criterion to determine whether a hydrodynamic description breaks down is to
find out the deviation of the eigenvalues in the corresponding fluid models
away from the reference eigenvalue which we denoted as $\nu_B$. We propose a criterion with the following two steps:
\begin{itemize}
\item {\bf Step 1}: If one cell is not in the kinetic regime, then if it is in the Euler regime, it will be classified as the NS regime if $|\nu_{B}-\nu_{Euler}|> \eta_0$.
\item {\bf Step 2}: For all the cells in the NS regime (with those from the Euler regime in Step 1), if $|\nu_{B}-\nu_{NS}|>\eta_1$, they will be classified as the kinetic regime.
\end{itemize}
The best choice of the thresholds $\eta_0$ and $\eta_1$ is problem dependent. In our numerical tests, we all take $\eta_0=10^{-2}$ and $\eta_1=10^{-1}$. In order to avoid classifying a smooth extremum (where $T_x$ and $u_x$ might both be zero) as the Euler regime (which may not be physically accurate), we propose to change a cell from the Euler regime into the NS regime if both of its neighbors are in the NS regime.

\begin{rem}
In our numerics, the first and second derivatives, such as $T_x$, $u_x$ and $u_{xx}$ when computing the eigenvalues, are approximated by a DG or LDG discretization with central fluxes.
\end{rem}

\subsubsection{From kinetic to fluid}
The criterion in this direction would simply be a comparison between the kinetic
density $f=M+\eps g$ and the truncated distribution $f_k$ whose moments match those of $f$ and $k$ is the order of the macroscopic model.
In particular, a kinetic description will be changed to a hydrodynamic closure of $k$ if
\beq
\|f(t,x,\cdot)-f_k(t,x,\cdot)\|_{L^2_M} \le \delta_0,
\eeq
with a weighted $L^2$ norm defined as $\| \cdot \|_{L^2_M}=\left(\int_{\mathbb{R}}|\cdot|^2/M_Udv/\rho\right)^{1/2}$.
For example, $f_1(t,x,v)=M$ and $\|f(t,x,v)-f_1(t,x,v)\|_{L^2_M}=\eps\|g(t,x,v)\|_{L^2_M}$,
while $f_2(t,x,v)=M-\eps A \frac{T_x}{\sqrt{T}} M$
and $\|f(t,x,v)-f_2(t,x,v)\|_{L^2_M}=\eps\|g(t,x,v)+A \frac{T_x}{\sqrt{T}} M\|_{L^2_M}$.
The criterion is implemented in the following two steps:
\begin{itemize}
\item {\bf Step 1}: For a cell not in the Euler regime (including both NS and kinetic), if $\eps\|g(t,x,v)\|_{L^2_M} < \delta_0$, it is added to the Euler regime.
\item {\bf Step 2}: Otherwise, if it is in the kinetic regime and $\eps\|g(t,x,v)+A \frac{T_x}{\sqrt{T}} M\|_{L^2_M} < \delta_0$, it is added to the NS regime.
\end{itemize}
In our numerical tests, we take $\delta_0=10^{-3}$.

\begin{rem}
In the hydrodynamic part including Euler and NS regimes, we only have the information of $U$, while in the kinetic regime both $U$ and $g$ are solved. In order to match the interfaces between the fluid and kinetic regimes, $g$ needs to be recovered in the fluid regime as the value at the thermal equilibrium $g=-A \frac{T_x}{\sqrt{T}} M$.
\end{rem}

\begin{rem}
Our hierarchical algorithm with regime indicators realizes an adaptive seamless coupling between hydrodynamic and kinetic solvers at different levels.
The scheme can be briefly sketched as follows:
\bit
\item
Initially at $T=0$ we start with all cells to be kinetic unless the initial conditions
are apparently to be fluid. For a high order Runge-Kutta method, the criteria are only
applied at the beginning of each time step.
\item 
At each intermediate stage from time level $t^n$ to $t^{n+1}$, we solve the Euler regime
with RK DG method, the NS regime with LDG method and the kinetic regime with the NDG-IMEX
method as described in Section \ref{sec:3.1} and Section \ref{sec:3.2}. {In all computational cells, macroscopic information ${U}$ are stored and being updated, while the microscopic components $g$ is stored and computed only in kinetic regimes. At boundary elements on hydrodynamic regimes (Euler or Navier-Stokes) border with the kinetic ones, we let $g=-A \frac{T_x}{\sqrt{T}} M$ as the boundary condition for microscopic component in kinetic regimes. }
\eit
{
We remark that the seamless coupling of different solvers are due to the compactness property of the DG method and the asymptotic equivalence of the kinetic solver with macroscopic ones. In particular, schemes on all three regimes are under the NDG-IMEX framework; the NDG-IMEX scheme is asymptotically equivalent to the DG and LDG scheme for Euler and Navier-Stokes systems. The macroscopic ones avoid computing several microscopic terms, leading to significant computational savings.}
\end{rem}

\section{Numerical Examples}
\label{sec:NDG-IMEX:4}

\setcounter{equation}{0}
\setcounter{figure}{0}
\setcounter{table}{0}

In this section, we will apply the regime indicators to the NDG-IMEX scheme for the micro-macro decomposed BGK equation \eqref{eq:mmdc} with constant $\eps$, and
\eqref{eq:mmd} with variable $\eps(x)$. We take the third order NDG-IMEX scheme with
$3$-point Gauss quadrature, corresponding to a Lagrangian polynomial basis of degree $2$.
The corresponding DG and LDG schemes in the compressible Euler and Navier-Stokes regimes
are also third order. The time step is chosen as $\Delta t= CFL h/\max(\Lambda, V_c)$, where $CFL=0.05$
and $\Lambda=\||u|+\sqrt{\gamma T}\|_{\infty}$ is the maximal absolute eigenvalue of $\partial F(U)/\partial U$ over the spatial domain.
The velocity domain $\Omega=[-V_c, V_c]$ is set to be large enough. The TVB limiter with the parameter $M_{tvb}=1$ is used and is only applied on $U_h$.

We will consider three different hierarchy schemes, where the hydrodynamic regime may
contain only the Euler regime, the NS regime, or both the Euler and NS regimes, which we
will denote as Euler-Kinetic, NS-Kinetic and Euler-NS-Kinetic.

\subsection{Sod problem}
First we consider the sod shock tube problem with initial conditions to be
\[
(\rho, u, T)=
\begin{cases}
(1,0,1) &\text{ if } x < 0.5 \\
(0.125,0,0.8) &\text{ if } x > 0.5.
\end{cases}
\]
on the domain $[-0.2, 1.2]$, $v\in[-4.5, 4.5]$.

We report the results for $\eps=10^{-2}$ in Figure \ref{fig11}
and for $\eps=10^{-3}$ in Figure \ref{fig21},
for Euler-kinetic, NS-kinetic and Euler-NS-kinetic at time $t=0.2$.
Here the hybrid schemes for all three cases are computed with $50$ cells, while the reference kinetic
solution and fluid solution are computed with $200$ cells. The TVB limiter is used with the parameter
to be $1$. It is observed that in regions where the hydrodynamic and kinetic solutions differ from each other, the kinetic solver is turned on and the hybrid solutions approximate the reference kinetic solutions well. The results for all three cases match each other very well. NS-kinetic has less kinetic cells
than Euler-kinetic. The Euler-NS-kinetic scheme with all three regimes together shows that NS regime well connects the Euler
and kinetic regimes.

For this problem, we also compare the computational cost for different indicators,
and also the full kinetic scheme. In Table \ref{tab1}, we can see that Euler-NS-kinetic can save $72\%$ for $\eps=10^{-2}$ and $87\%$ for $\eps=10^{-3}$ as compared to the full kinetic scheme, while NS-kinetic saves $74\%$ for $\eps=10^{-2}$ and $85\%$ for $\eps=10^{-3}$. Euler-kinetic is a little higher, but still $55\%$ for $\eps=10^{-2}$ and $70\%$ for $\eps=10^{-3}$. This is due to that the kinetic solution is very close to the Euler solution, especially when $\eps=10^{-3}$, so the regime indicators are very efficient for this problem. The Euler-NS-kinetic performs almost the same as the NS-kinetic one, costing a little higher for $\eps=10^{-2}$ and a little lower for $\eps=10^{-3}$.

\begin{table}[!h]
\label{tab1}
\begin{center}
\caption{{Comparison of the computational time (seconds) for different indicators. Sod problem. }
}
\bigskip
\begin{tabular}{|c|c|c|c|c|}
\hline
\cline{1-5} method  & Euler-NS-kinetic & NS-kinetic & Euler-kinetic & Full kinetic  \\
\hline
\cline{1-5}  $\eps=10^{-2}$ & 8.38 & 7.81 & 13.40 & 30.19 \\
\hline
\cline{1-5}  $\eps=10^{-3}$ & 3.78 & 4.28 & 8.28 & 28.57 \\
\hline
\end{tabular}
\end{center}
\end{table}

\begin{figure}[!h]
\centering
\includegraphics[totalheight=1.6in, width = 2.1in]{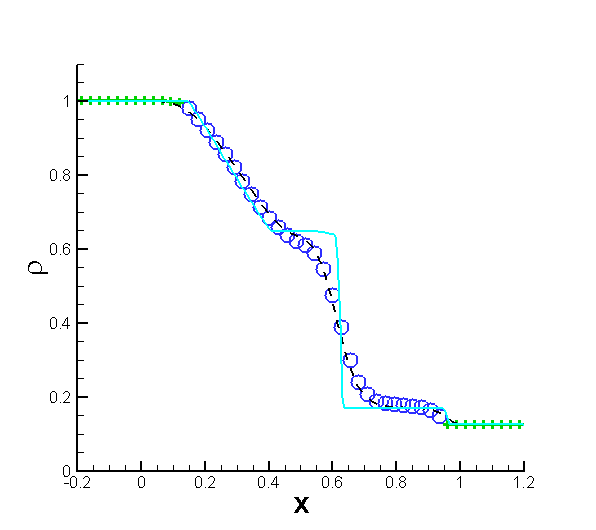},
\includegraphics[totalheight=1.6in, width = 2.1in]{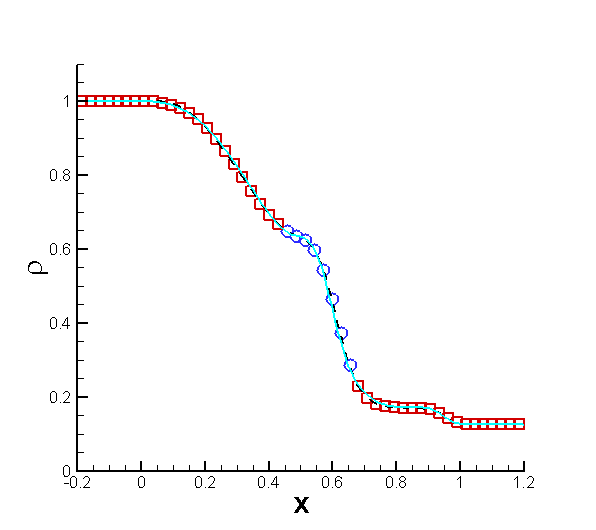},
\includegraphics[totalheight=1.6in, width = 2.1in]{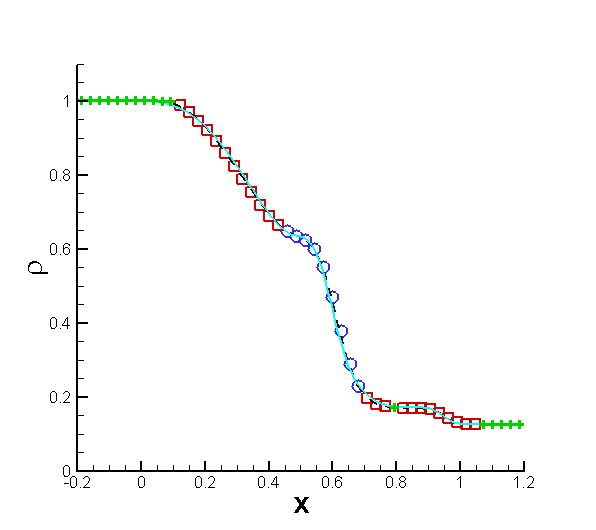} \\
\includegraphics[totalheight=1.6in, width = 2.1in]{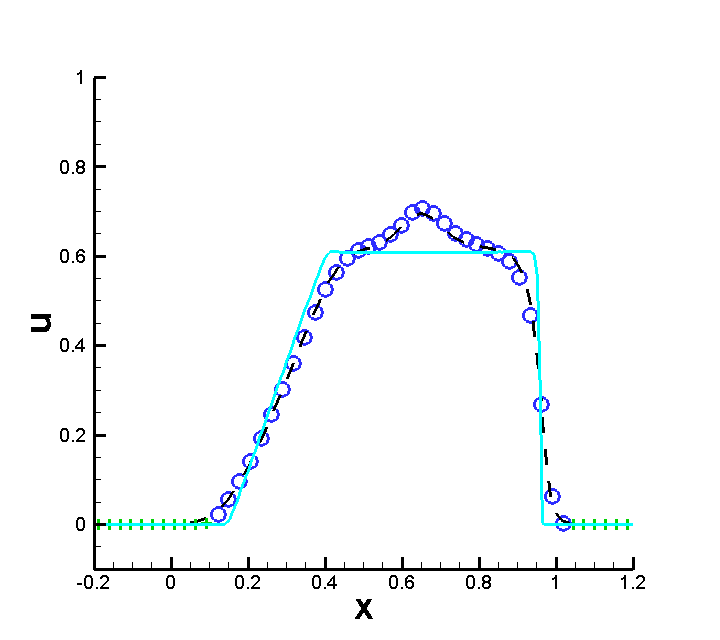},
\includegraphics[totalheight=1.6in, width = 2.1in]{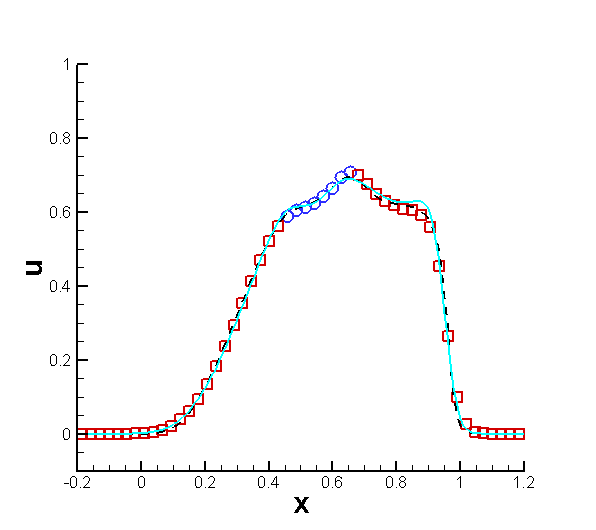},
\includegraphics[totalheight=1.6in, width = 2.1in]{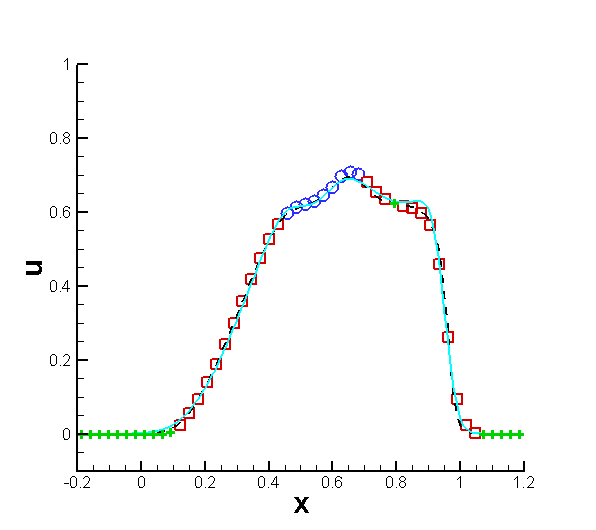} \\
\includegraphics[totalheight=1.6in, width = 2.1in]{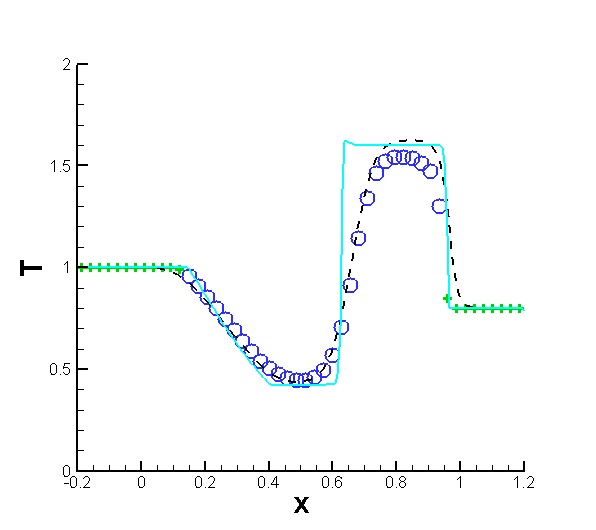},
\includegraphics[totalheight=1.6in, width = 2.1in]{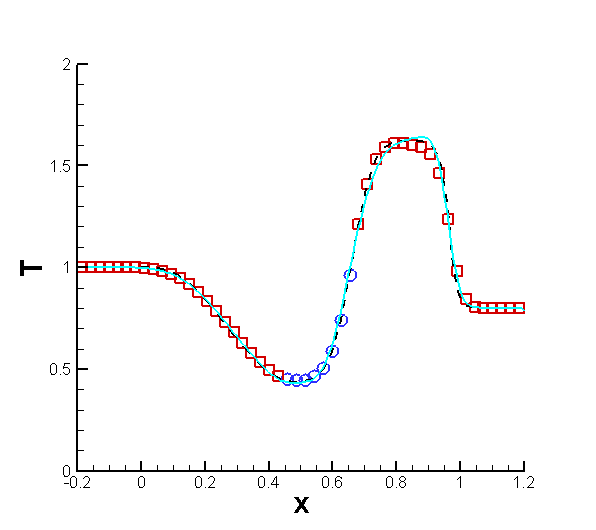},
\includegraphics[totalheight=1.6in, width = 2.1in]{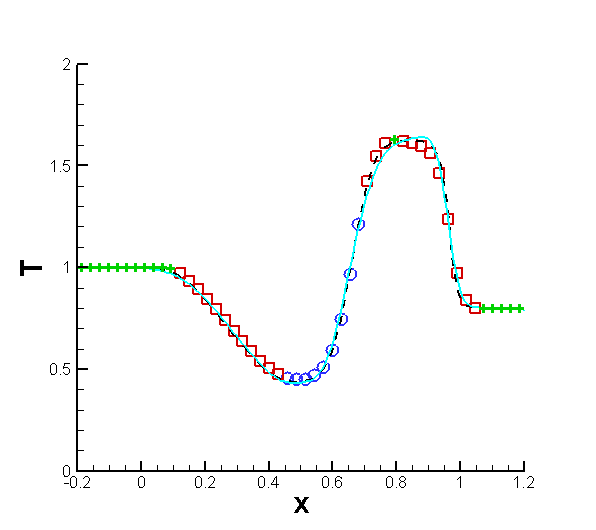},\\
\includegraphics[totalheight=1.6in, width = 2.1in]{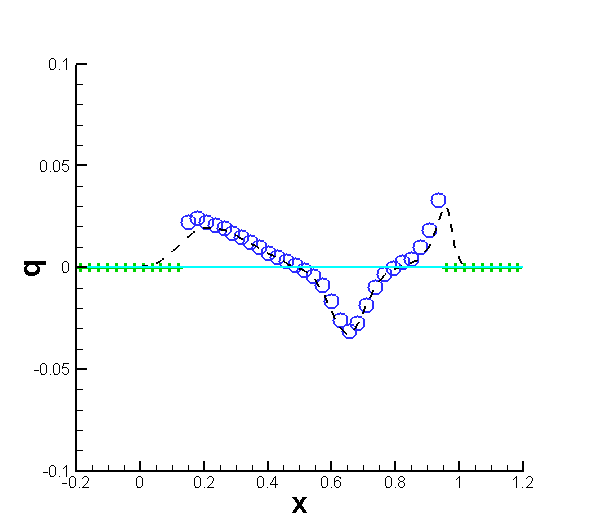},
\includegraphics[totalheight=1.6in, width = 2.1in]{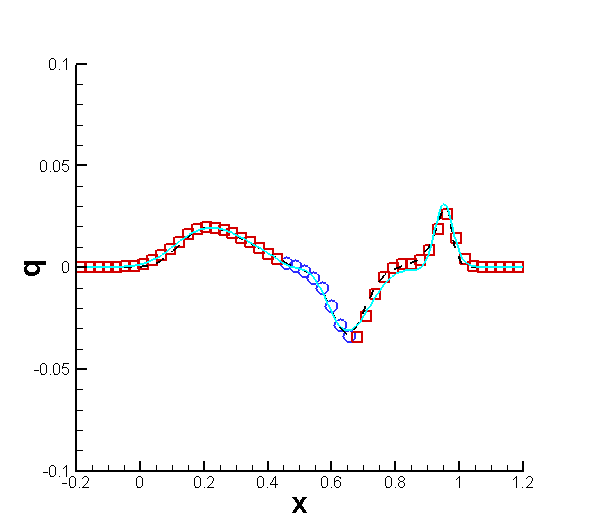},
\includegraphics[totalheight=1.6in, width = 2.1in]{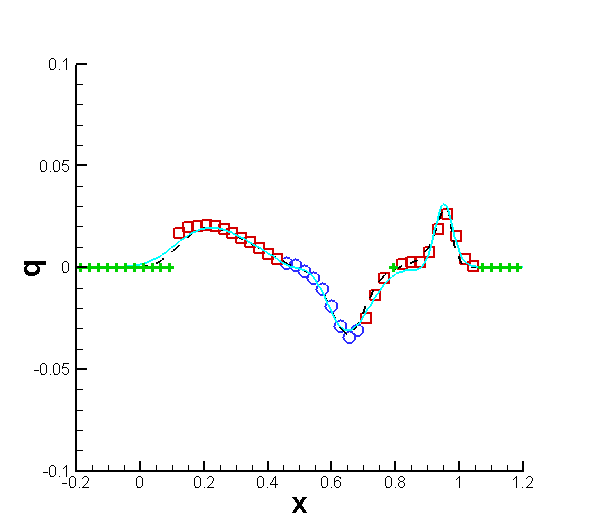}
\caption{Sod problem at time $t=0.2$ on the domain $[-0.2,1.2]\times[-4.5,4.5]$.
Symbols: $N_x=50$ and $N_v=100$ with NDG3. From left to right: Euler-kinetic, NS-kinetic and Euler-NS-kinetic. From top to bottom,
the density $\rho$, the mean velocity $u$, the temperature $T$ and the heat flux $q$.
Symbols: `+' is Euler, circle is kinetic. Solid line: the Euler reference solution for Euler-kinetic and the NS reference solution for the other two. Dashed line: the kinetic reference solution. $\eps=10^{-2}$.}
\label{fig11}
\end{figure}

\begin{figure}[!h]
\centering
\includegraphics[totalheight=1.6in, width = 2.1in]{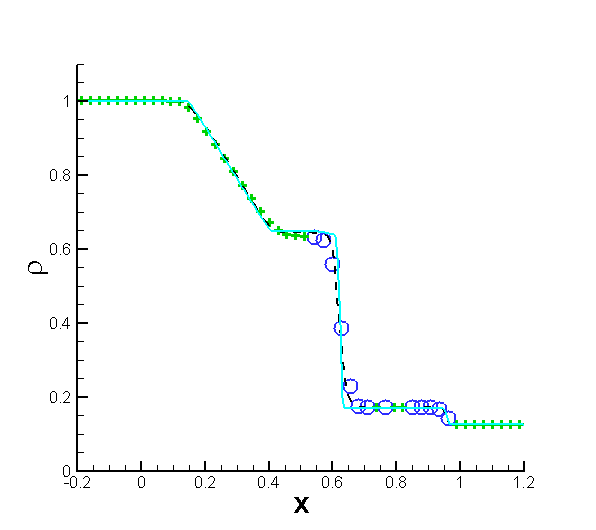},
\includegraphics[totalheight=1.6in, width = 2.1in]{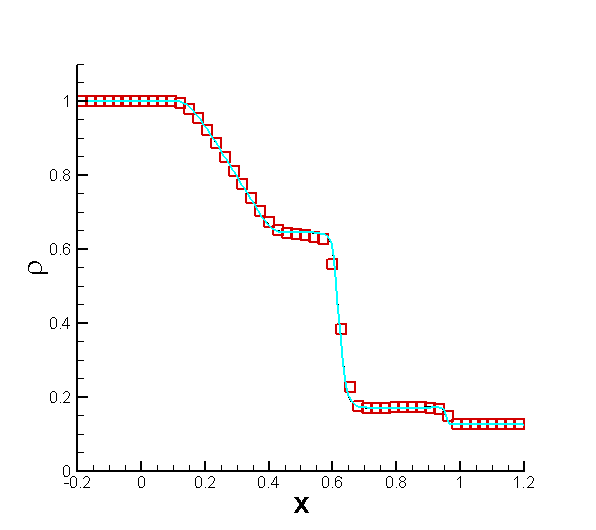},
\includegraphics[totalheight=1.6in, width = 2.1in]{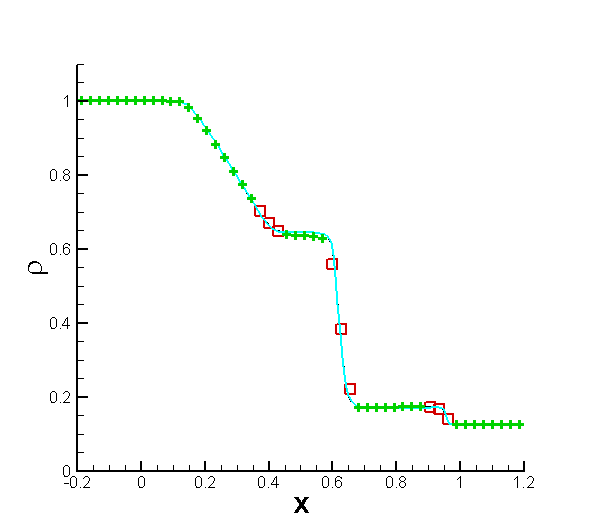} \\
\includegraphics[totalheight=1.6in, width = 2.1in]{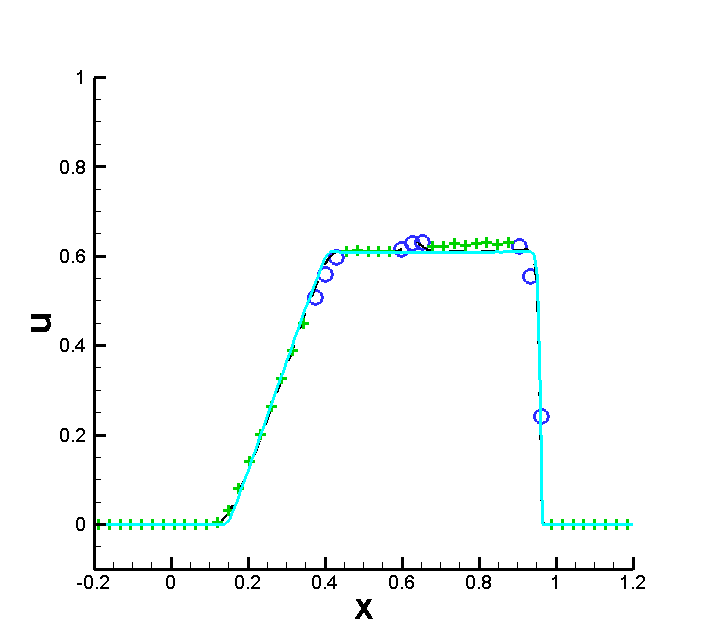},
\includegraphics[totalheight=1.6in, width = 2.1in]{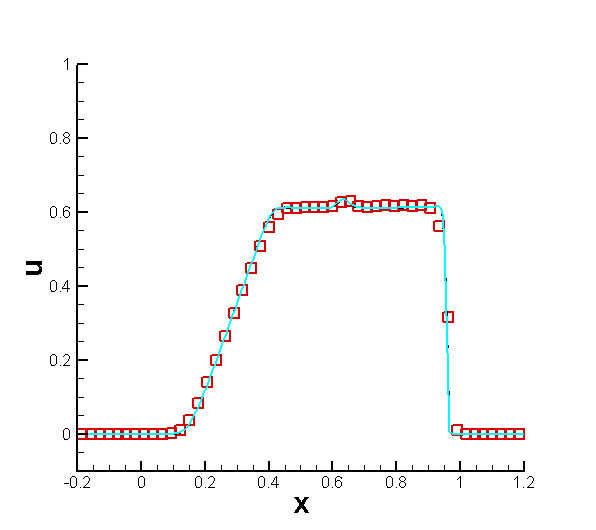},
\includegraphics[totalheight=1.6in, width = 2.1in]{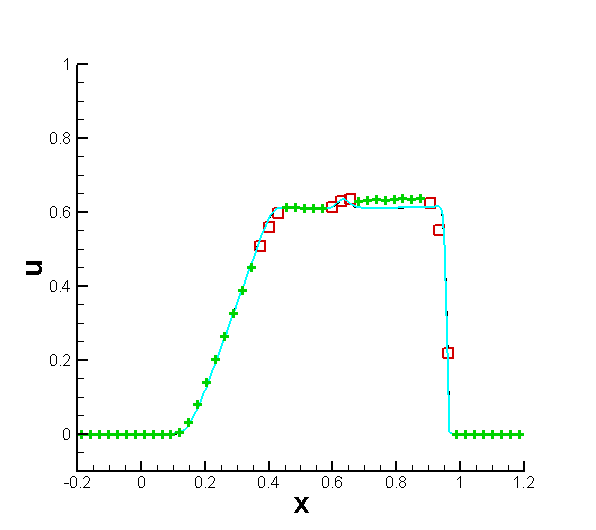} \\
\includegraphics[totalheight=1.6in, width = 2.1in]{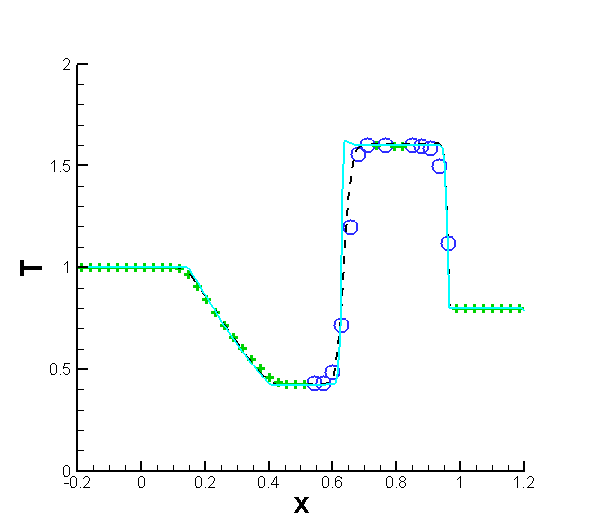},
\includegraphics[totalheight=1.6in, width = 2.1in]{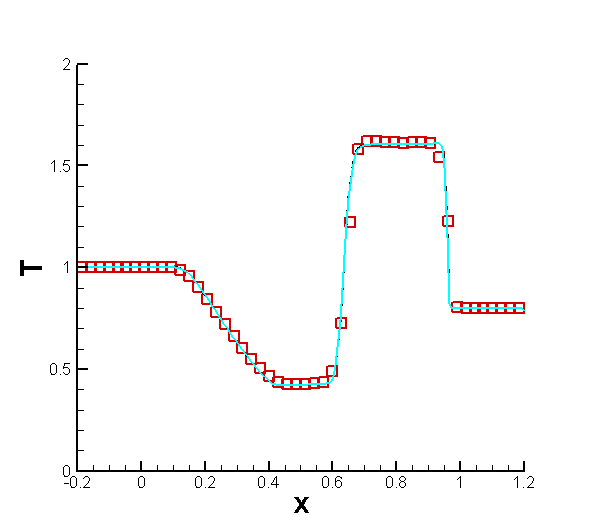},
\includegraphics[totalheight=1.6in, width = 2.1in]{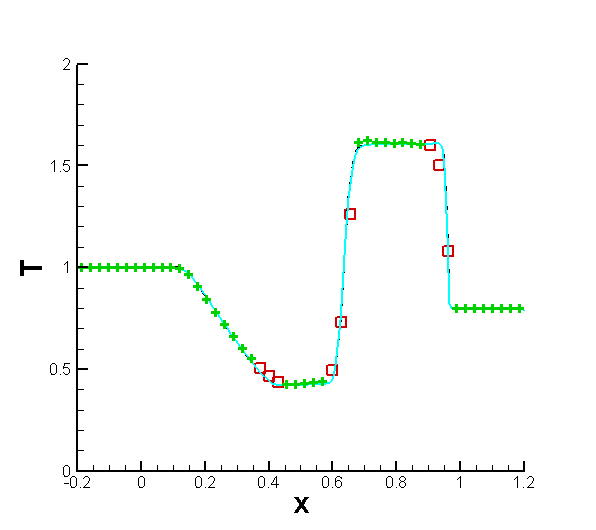},\\
\includegraphics[totalheight=1.6in, width = 2.1in]{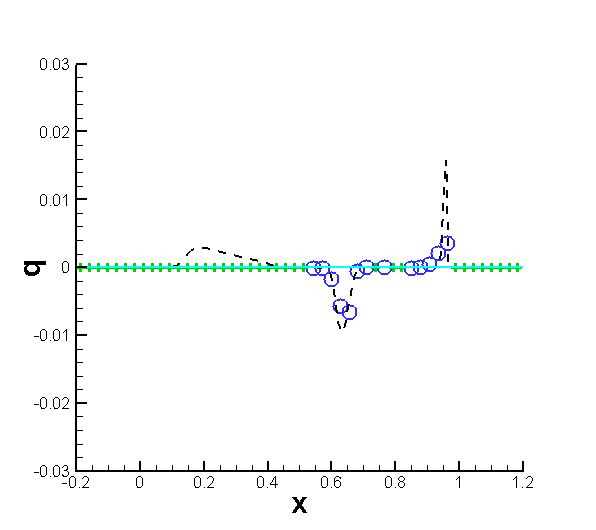},
\includegraphics[totalheight=1.6in, width = 2.1in]{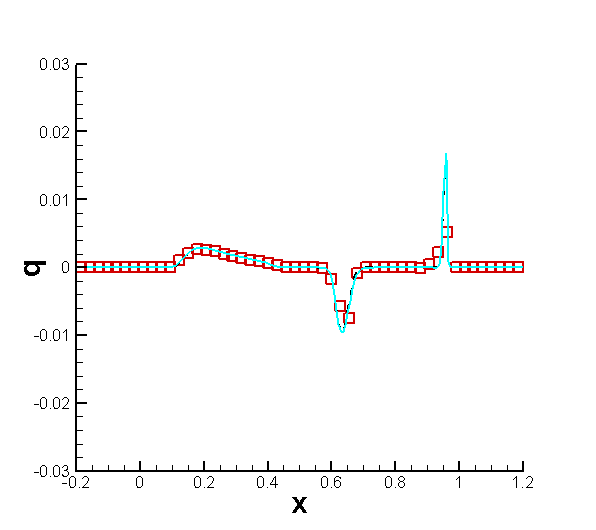},
\includegraphics[totalheight=1.6in, width = 2.1in]{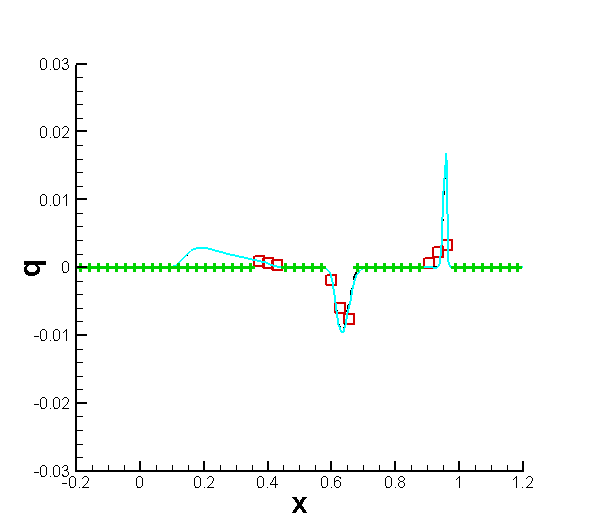}
\caption{Sod problem at time $t=0.2$ on the domain $[-0.2,1.2]\times[-4.5,4.5]$.
Symbols: $N_x=50$ and $N_v=100$ with NDG3. From left to right: Euler-kinetic, NS-kinetic and Euler-NS-kinetic. From top to bottom,
the density $\rho$, the mean velocity $u$, the temperature $T$ and the heat flux $q$.
Symbols: `+' is Euler, circle is kinetic. Solid line: the Euler reference solution for Euler-kinetic and the NS reference solution for the other two. Dashed line: the kinetic reference solution. $\eps=10^{-3}$.}
\label{fig21}
\end{figure}

\subsection{Blast wave problem}
For the blast wave problem, the initial condition is given by
\[
(\rho, u, T)=
\begin{cases}
(1,1,2) &\text{ if } x < 0.2 \\
(1,0,0.25) &\text{ if } 0.2\le x \le 0.8 \\
(1,-1,2) &\text{ if } x > 0.8.
\end{cases}
\]
with reflective boundary condition in the $x$ direction on the domain $[0,1]\times[-9,9]$.

Similarly the hybrid scheme is computed with $50$ cells, while the reference kinetic
solution and fluid solution are computed with $200$ cells. The TVB limiter is used with the parameter to be $1$. Here we report the results for $\eps=10^{-2}$ and $\eps=10^{-3}$ for Euler-NS-kinetic in Figures \ref{fig33} and \ref{fig43}, respectively.


For this problem, when $\eps=10^{-2}$, the kinetic solution deviate slightly away from the fluid solution, and most computational cells are assigned into the kinetic regime. However, when $\eps$ becomes smaller, that is when $\eps=10^{-3}$, the kinetic solution is getting close to the hydrodynamic solution, it can be observed from Figure \ref{fig43} that the Euler and NS solvers are turned on in larger regions, leading to computational savings. If one compares Figure \ref{fig43}, ``best" solvers that well balance computational effectiveness (in capturing reference solutions) and efficiency (in saving computational time) are adaptively chosen by the criteria.

Similar to the Sod problem, we compare the computational cost for the blast wave problem
in Table \ref{tab2} for three different hierarchy schemes. As we can see, when $\eps=10^{-2}$, since the solution is mostly in the kinetic regime. Euler-NS-kinetic and NS-kinetic only save around $10\%\sim15\%$, while Euler-kinetic takes even more computational time than the full kinetic scheme. Also the computational time of Euler-NS-kinetic is slightly more than that of the NS-kinetic. This is due to the fact that computing the regime indicators takes extra computational time. When $\eps=10^{-3}$, the solution becomes closer to the NS solution, Euler-NS-kinetic and NS-kinetic can save up to $80\%$, while Euler-kinetic can only save around $60\%$. Also Euler-NS-kinetic takes more time than NS-kinetic. For the 1D problem, the compressible Navier-Stokes equations \eqref{cns-bgk-1d} has only one extra term than the compressible Euler equations.
We might have taken slightly more time to compute the regime indicators and the logic
decisions for the Euler-NS-kinetic approach, while computing the extra term does not take much time for 1D cases. We would expect that the savings of the Euler-NS-kinetic would become significant for high dimensional problems.

\begin{table}[!h]
\label{tab2}
\begin{center}
\caption{{Comparison of the computational cost (seconds) for different indicators. Blast wave problem. }
}
\bigskip
\begin{tabular}{|c|c|c|c|c|}
\hline
\cline{1-5} method  & Euler-NS-kinetic & NS-kinetic & Euler-kinetic & Full kinetic  \\
\hline
\cline{1-5}  $\eps=10^{-2}$ & 85.74 & 80.98 & 102.42 & 96.72 \\
\hline
\cline{1-5}  $\eps=10^{-3}$ & 17.22 & 16.77 & 37.97  & 98.72 \\
\hline
\end{tabular}
\end{center}
\end{table}

\begin{figure}[!h]
\centering
\includegraphics[totalheight=1.6in, width = 2.1in]{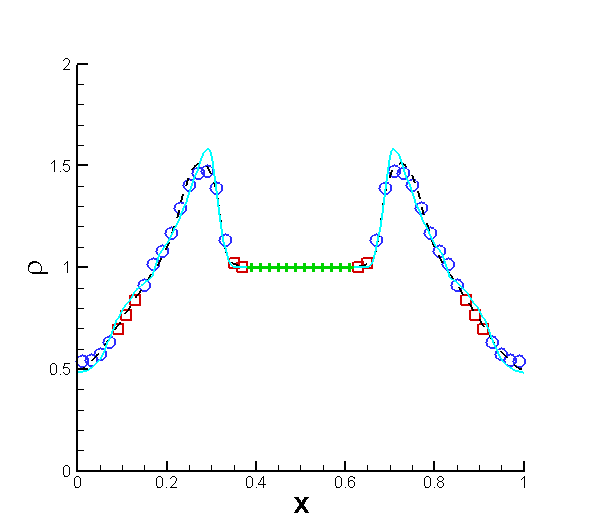},
\includegraphics[totalheight=1.6in, width = 2.1in]{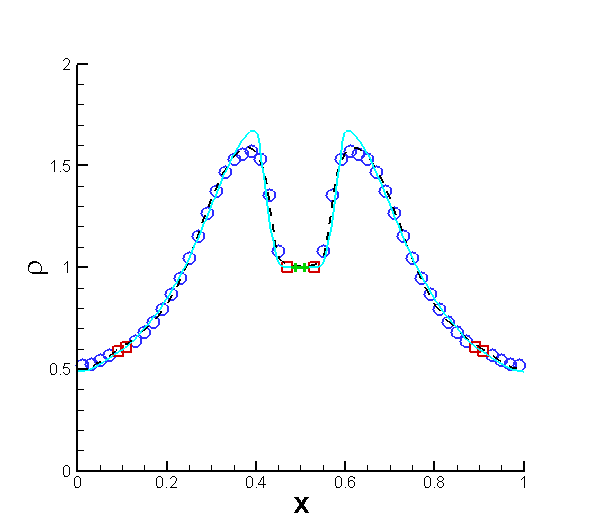},
\includegraphics[totalheight=1.6in, width = 2.1in]{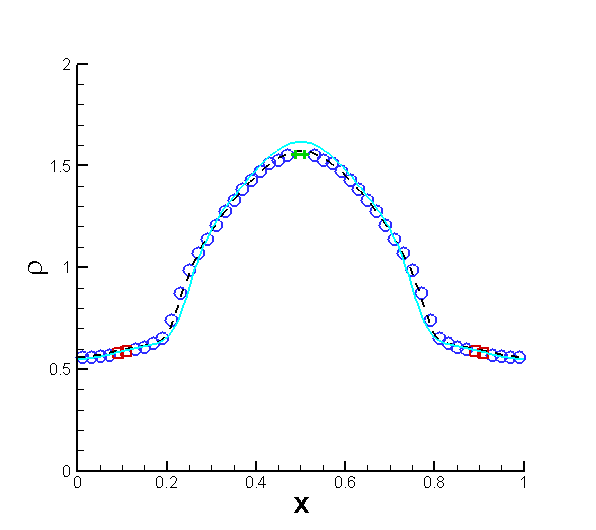} \\
\includegraphics[totalheight=1.6in, width = 2.1in]{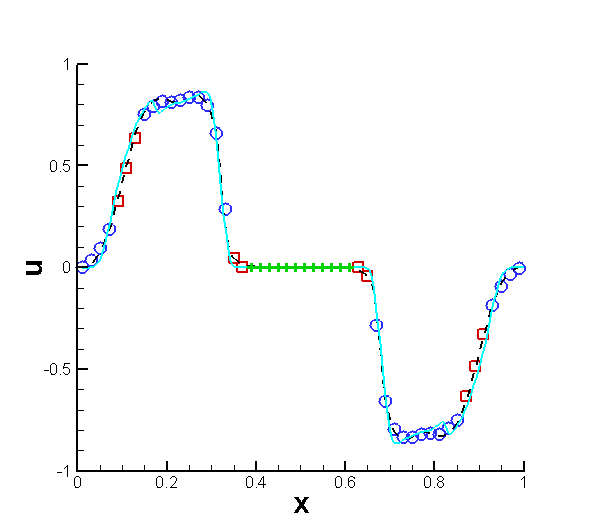},
\includegraphics[totalheight=1.6in, width = 2.1in]{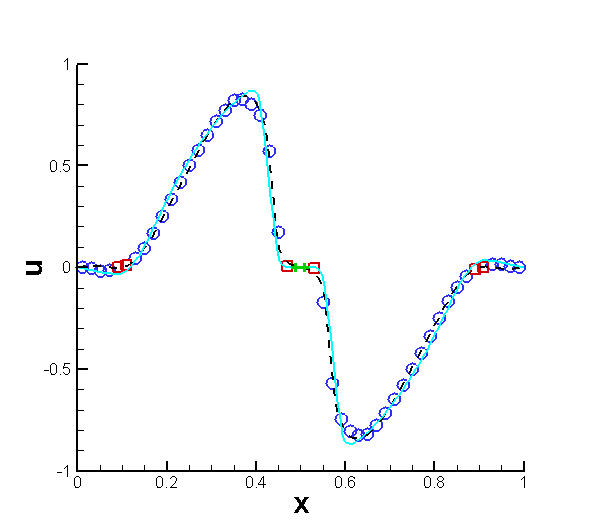},
\includegraphics[totalheight=1.6in, width = 2.1in]{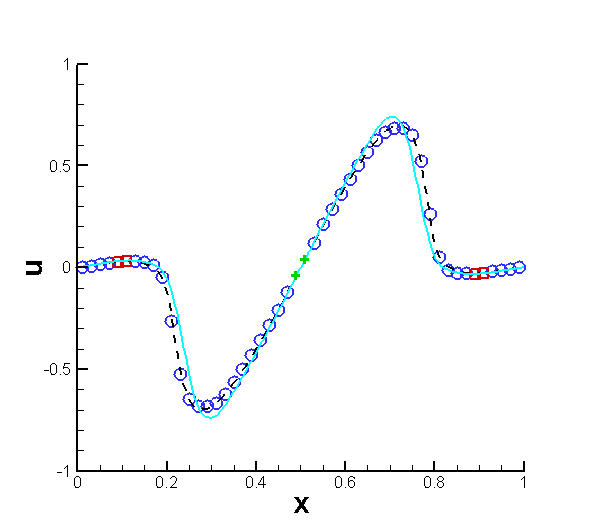} \\
\includegraphics[totalheight=1.6in, width = 2.1in]{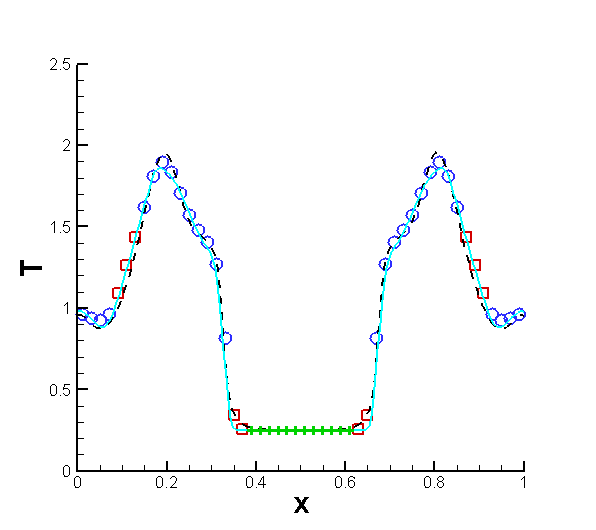},
\includegraphics[totalheight=1.6in, width = 2.1in]{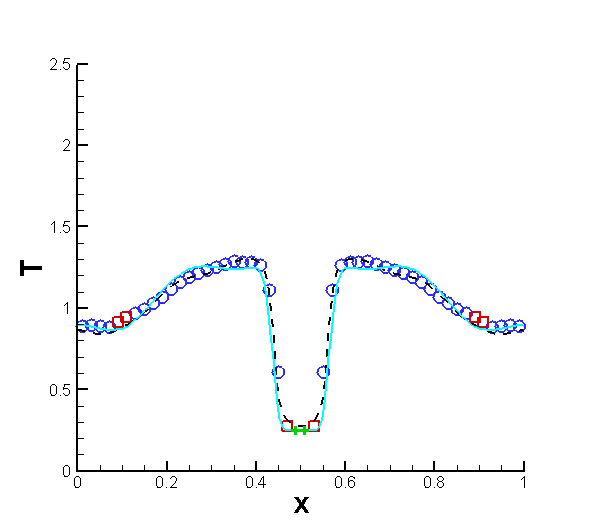},
\includegraphics[totalheight=1.6in, width = 2.1in]{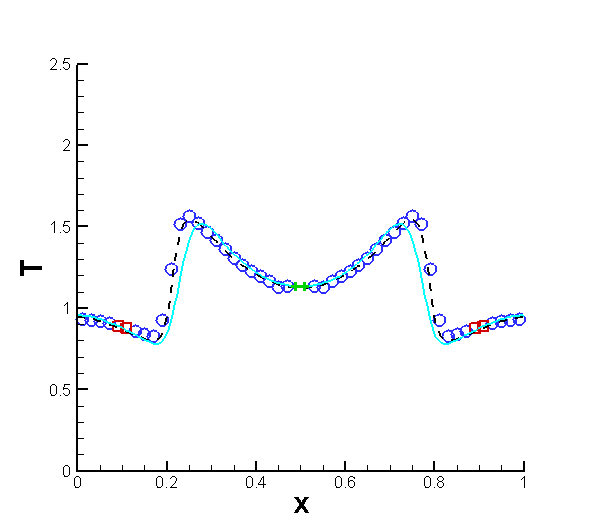},\\
\includegraphics[totalheight=1.6in, width = 2.1in]{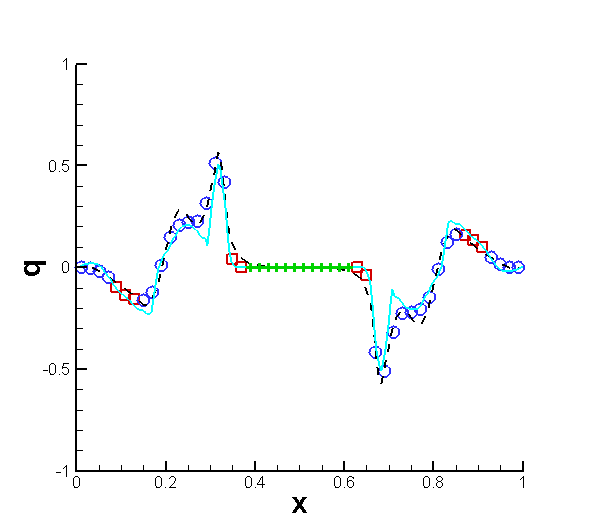},
\includegraphics[totalheight=1.6in, width = 2.1in]{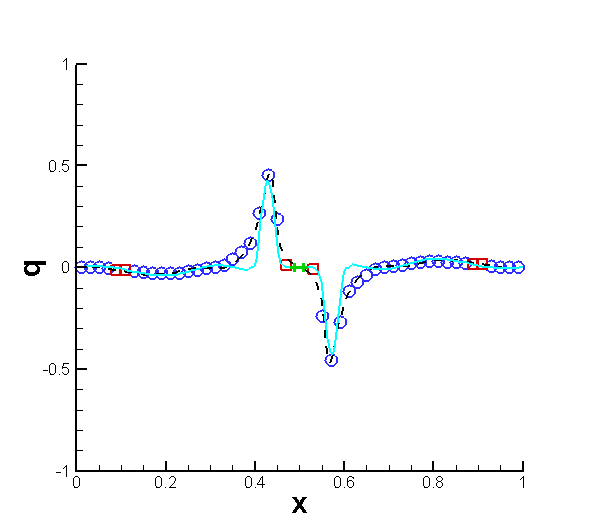},
\includegraphics[totalheight=1.6in, width = 2.1in]{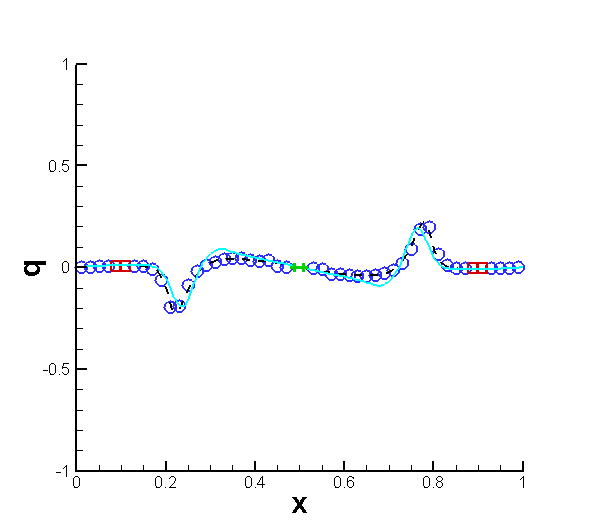}
\caption{Blast wave problem on the domain $[0,1]\times[-9,9]$.
Symbols: $N_x=50$ and $N_v=100$ with NDG3. From left to right: time
$t=0.05, 0.1, 0.25$. From top to bottom,
the density $\rho$, the mean velocity $u$, the temperature $T$ and the heat flux $q$.
Symbols: `+' is Euler, square is NS, circle is kinetic. Solid line: the NS reference solution. Dashed line: the kinetic reference solution. $\eps=10^{-2}$, Euler-NS-kinetic.}
\label{fig33}
\end{figure}

\begin{figure}[!h]
\centering
\includegraphics[totalheight=1.6in, width = 2.1in]{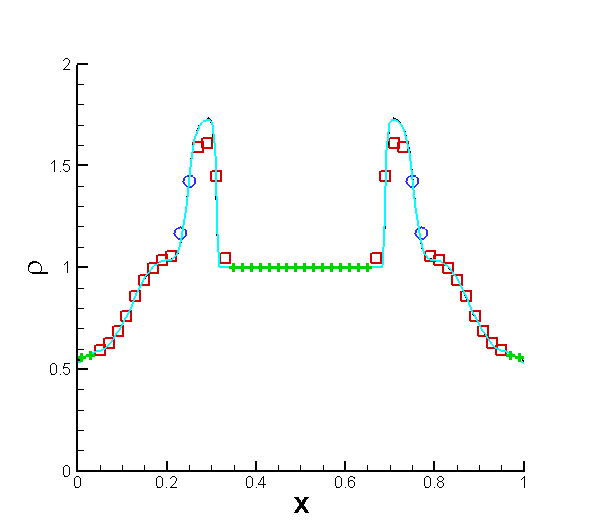},
\includegraphics[totalheight=1.6in, width = 2.1in]{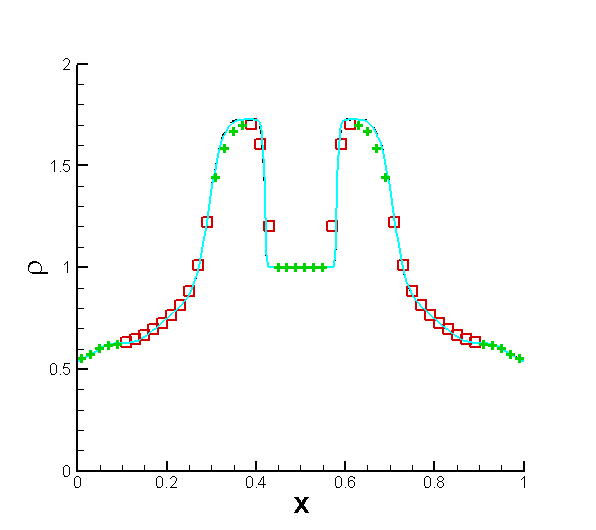},
\includegraphics[totalheight=1.6in, width = 2.1in]{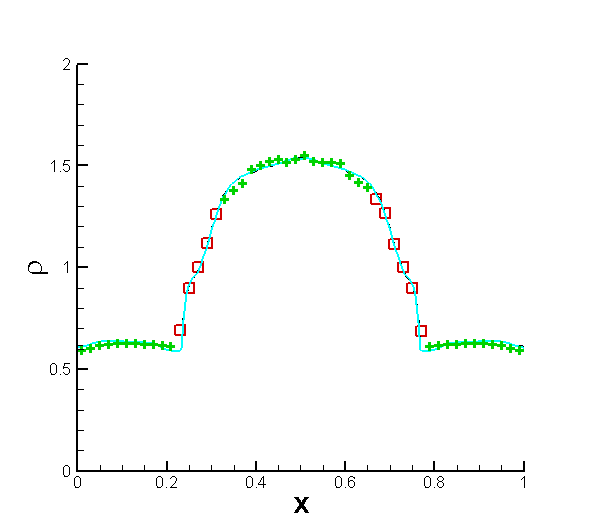} \\
\includegraphics[totalheight=1.6in, width = 2.1in]{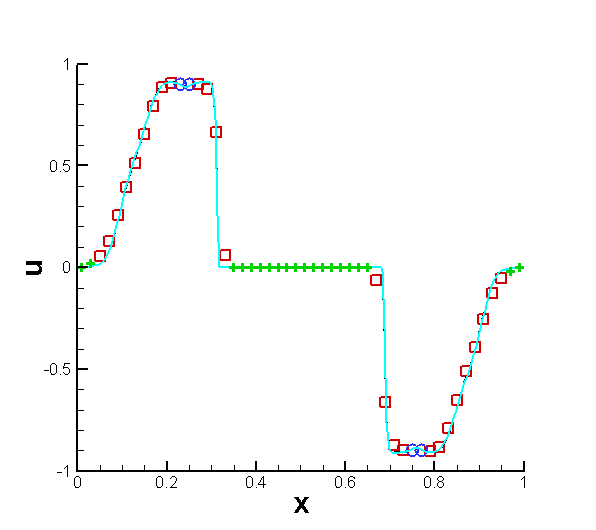},
\includegraphics[totalheight=1.6in, width = 2.1in]{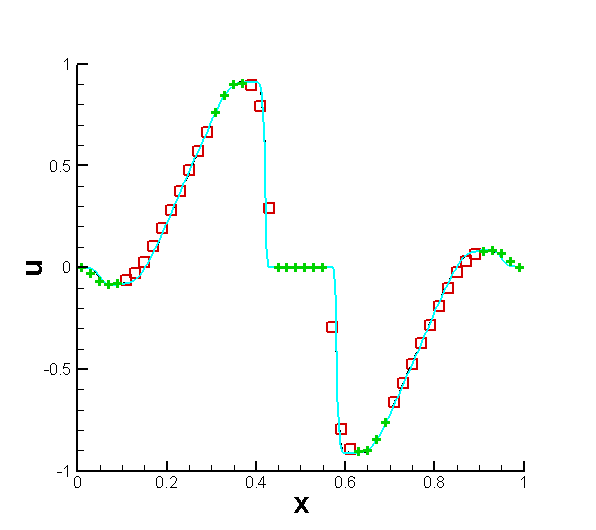},
\includegraphics[totalheight=1.6in, width = 2.1in]{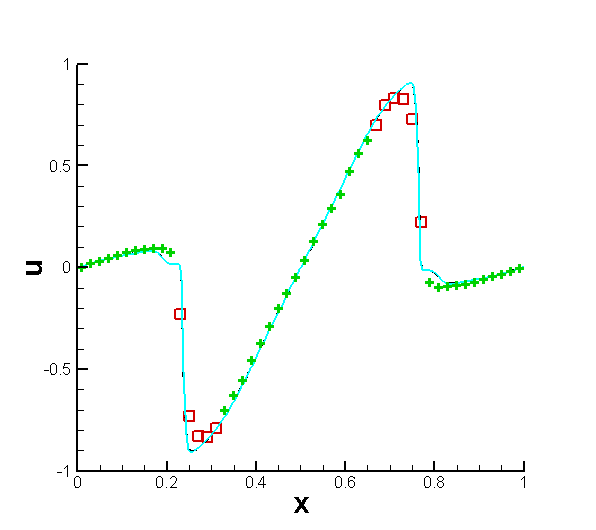} \\
\includegraphics[totalheight=1.6in, width = 2.1in]{./pic_hybrid/blu21.png},
\includegraphics[totalheight=1.6in, width = 2.1in]{./pic_hybrid/blu22.png},
\includegraphics[totalheight=1.6in, width = 2.1in]{./pic_hybrid/blu23.png},\\
\includegraphics[totalheight=1.6in, width = 2.1in]{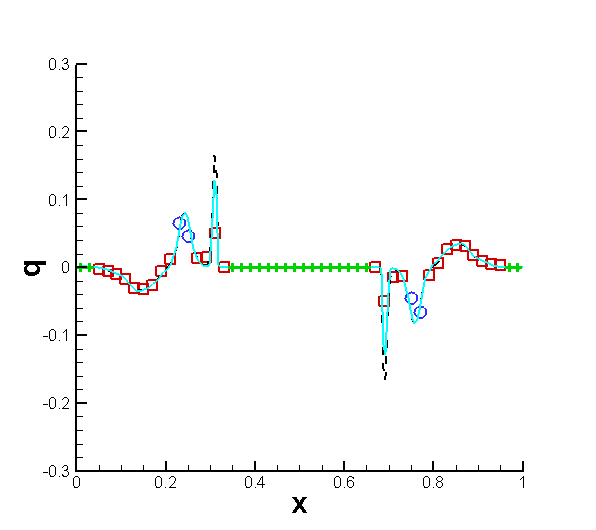},
\includegraphics[totalheight=1.6in, width = 2.1in]{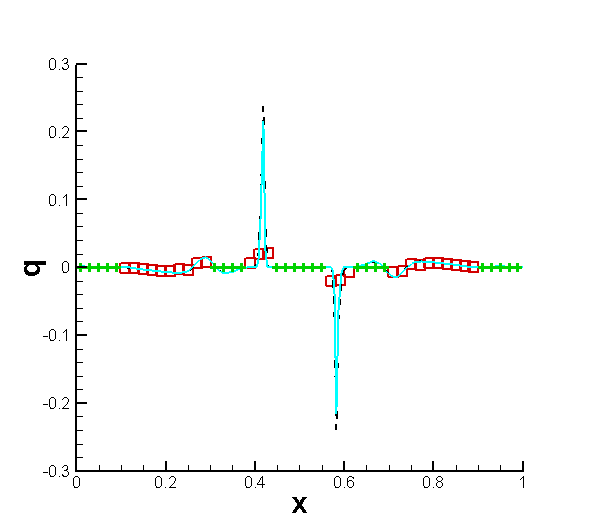},
\includegraphics[totalheight=1.6in, width = 2.1in]{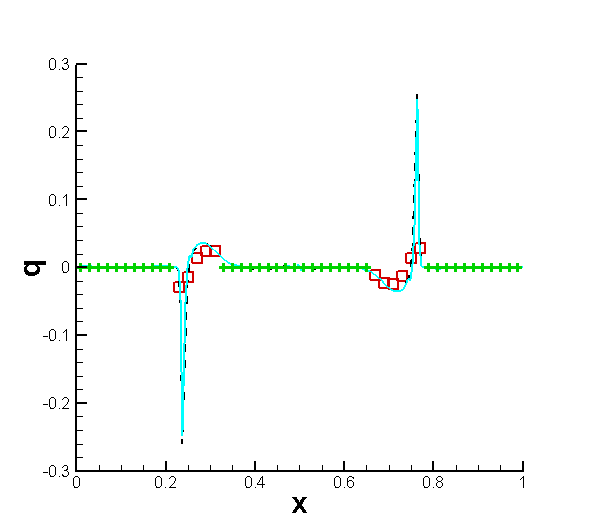}
\caption{Blast wave problem on the domain $[0,1]\times[-9,9]$.
Symbols: $N_x=50$ and $N_v=100$ with NDG3. From left to right: time
$t=0.05, 0.1, 0.25$. From top to bottom,
the density $\rho$, the mean velocity $u$, the temperature $T$ and the heat flux $q$.
Symbols: `+' is Euler, square is NS, circle is kinetic. Solid line: the NS reference solution. Dashed line: the kinetic reference solution. $\eps=10^{-3}$, Euler-NS-kinetic.}
\label{fig43}
\end{figure}


\subsection{Mixed regime problem}

We consider an example with a variable $\eps(x)$,
\begin{equation}
\eps(x)=\eps_0+\frac12\Big(\tanh(1- a_0 x)+\tanh(1+a_0 x)\Big),
\label{epsx}
\end{equation}
with $\eps_0=10^{-3}$. In the middle area of $x\in[-0.18, 0.18]$, since $\eps=\mathcal{O}(1)$, this part is always to be in the kinetic regime.
The initial distribution function $f$ is far away from the Maxwellian, which is
\begin{eqnarray}
f(x,v,0)=\frac{\tilde \rho}{2( 2 \pi \tilde T)^{1/2}}\left[\exp\left(-\frac{|v-\tilde u|^2}{2 \tilde T}\right)
+\exp\left(-\frac{|v+\tilde u|^2}{2\tilde T}\right)\right],
\end{eqnarray}
with
\begin{equation}
\tilde \rho(x)=1+0.875 \sin(w x), \quad \tilde T(x)=0.5+0.4 \sin(\omega x), \quad \tilde u(x)=0.75,
\end{equation}
on the spatial domain $x\in [-L, L]$, where $\omega=\pi/L$ and $L=0.5$. The initial macroscopic variables are
\begin{equation}
\rho(x,0)=\tilde \rho(x), \quad u(x,0)=0, \quad T(x,0)=\tilde T(x)+ \tilde u(x)^2,
\end{equation}
and the initial Maxwellian distribution is
\begin{equation}
M_U(x,v,0)=\frac{\rho(x,0)}{(2\pi T(x,0))^{1/2}}\exp\left(-\frac{|v-u(x,0)|^2}{2T(x,0)}\right).
\end{equation}
Periodic boundary conditions are used for both $U$ and $g$ in the $x$ direction. The velocity domain is taken to be $\Omega_v=[-10, 10]$.

We report the results with $50$ cells for the computational solution and $200$ cells
for the reference solutions in Figures \ref{fig53}. Here the NS reference
solution is obtained under the time step $\Delta t = \mathcal{O}(h^2)$. However, in the
middle region, since $\eps(x)=\mathcal{O}(1)$, the NS solution deviates far away from
the kinetic solution and the heat flux oscillates greatly as seen in Figure \ref{fig53}.
The adaptive algorithm can capture the kinetic solution well. For the computational cost,
from Table \ref{tab3}, we can observe similar result as the blast wave problem for $\eps=10^{-3}$. $40\%\sim 45\%$ of the computational cost for the Euler-NS-kinetic and
NS-kinetic methods, but only $20\%$ for the Euler-kinetic one. 

\begin{table}[!h]
\label{tab3}
\begin{center}
\caption{{Comparison of the computational cost ($seconds$) for different indicators. Mixed regime problem. }
}
\bigskip
\begin{tabular}{|c|c|c|c|c|}
\hline
\cline{1-5} method  & Euler-NS-kinetic & NS-kinetic & Euler-kinetic & Full kinetic  \\
\hline
\cline{1-5}  $\eps_0=10^{-3}$ & 362.52 & 343.44 & 512.68 & 631.07 \\
\hline
\end{tabular}
\end{center}
\end{table}

\begin{figure}[!h]
\centering
\includegraphics[totalheight=1.6in, width = 2.1in]{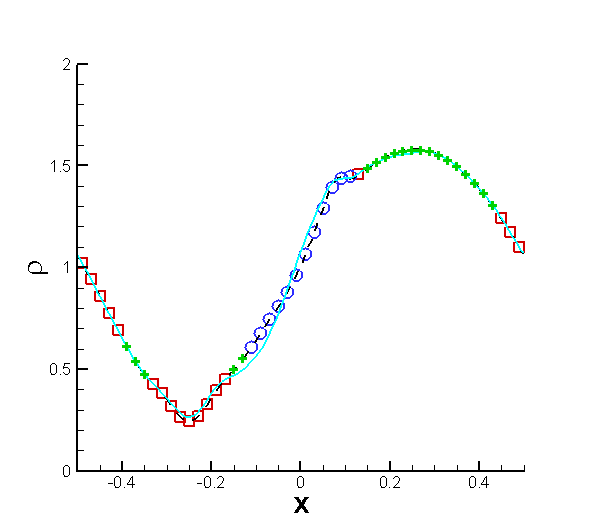},
\includegraphics[totalheight=1.6in, width = 2.1in]{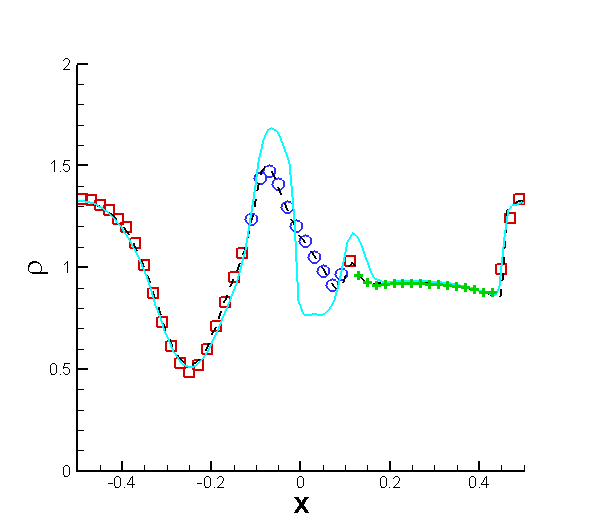},
\includegraphics[totalheight=1.6in, width = 2.1in]{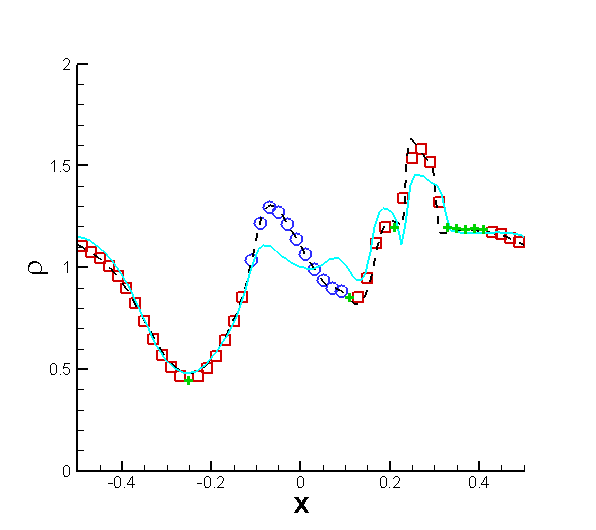} \\
\includegraphics[totalheight=1.6in, width = 2.1in]{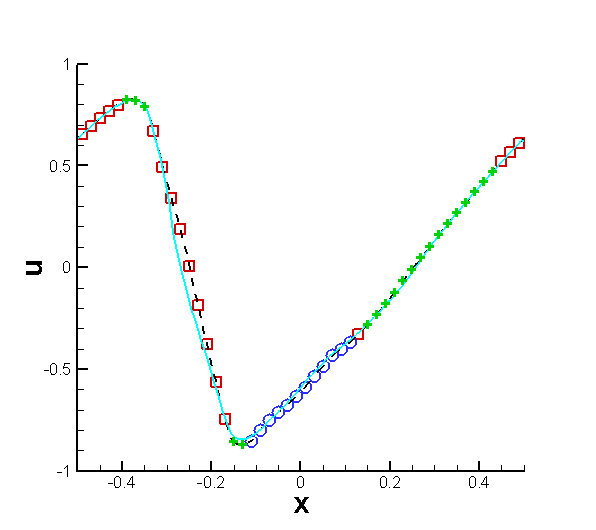},
\includegraphics[totalheight=1.6in, width = 2.1in]{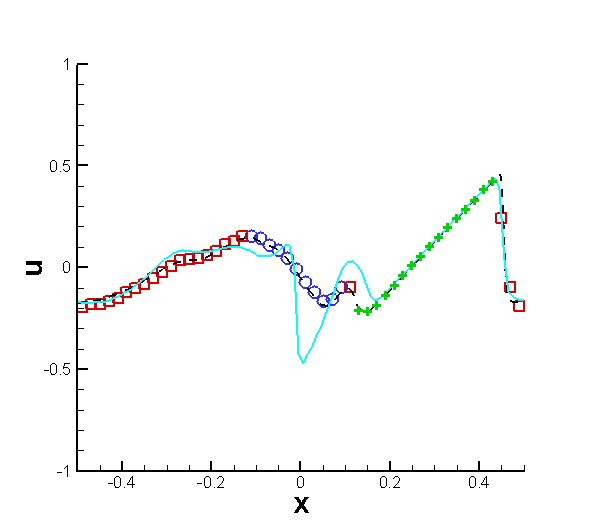},
\includegraphics[totalheight=1.6in, width = 2.1in]{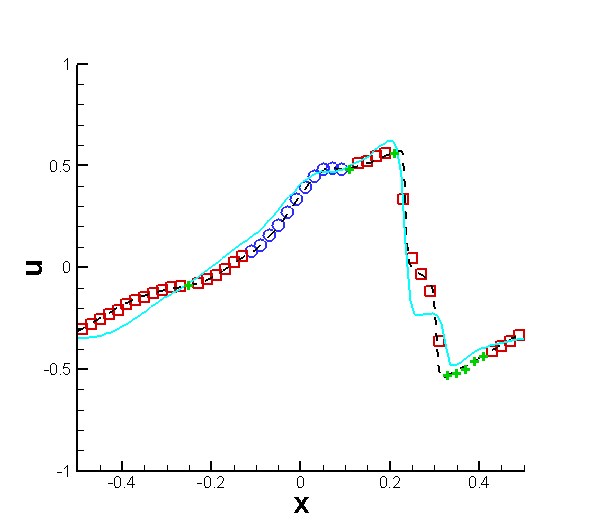} \\
\includegraphics[totalheight=1.6in, width = 2.1in]{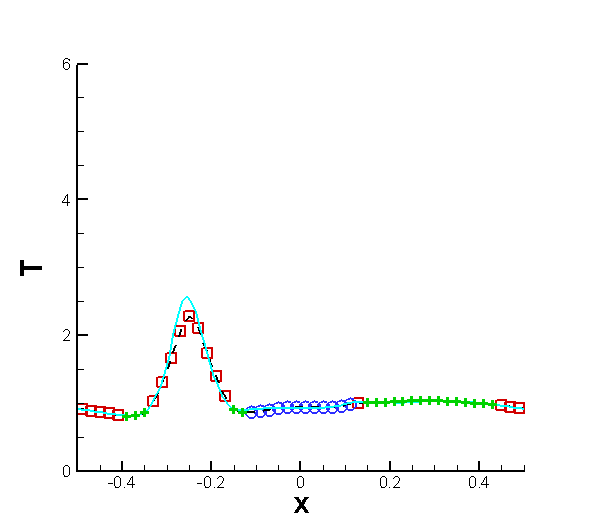},
\includegraphics[totalheight=1.6in, width = 2.1in]{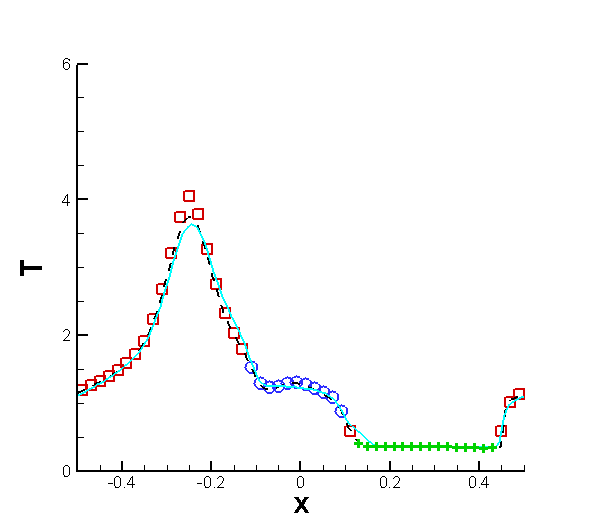},
\includegraphics[totalheight=1.6in, width = 2.1in]{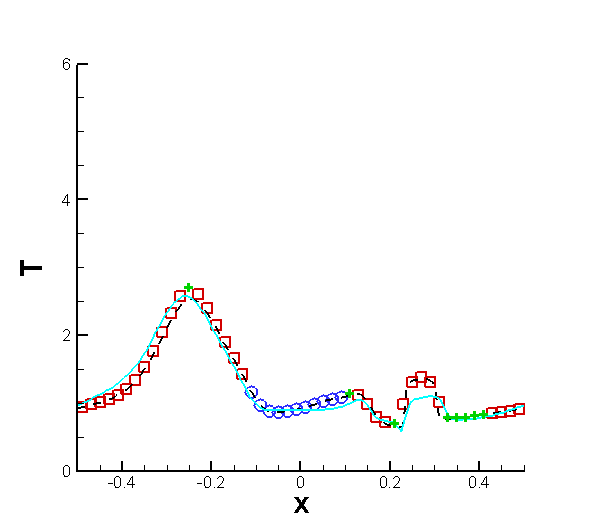},\\
\includegraphics[totalheight=1.6in, width = 2.1in]{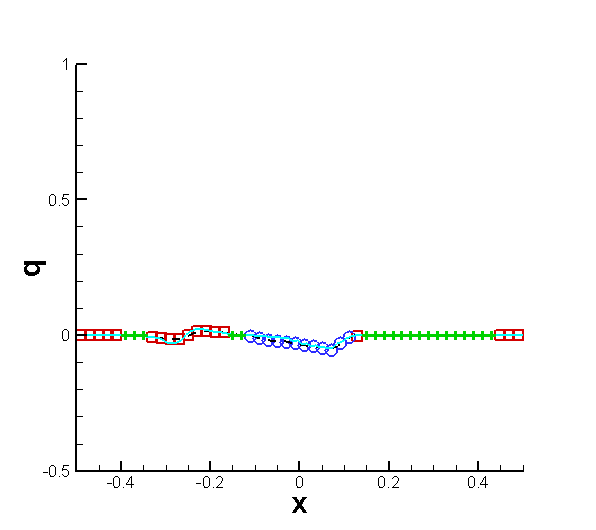},
\includegraphics[totalheight=1.6in, width = 2.1in]{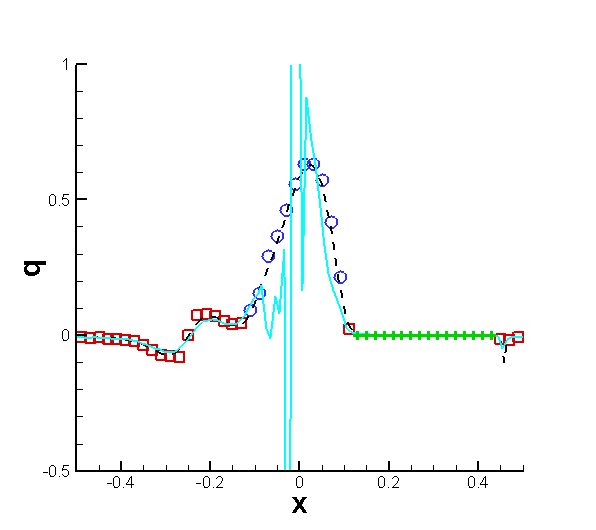},
\includegraphics[totalheight=1.6in, width = 2.1in]{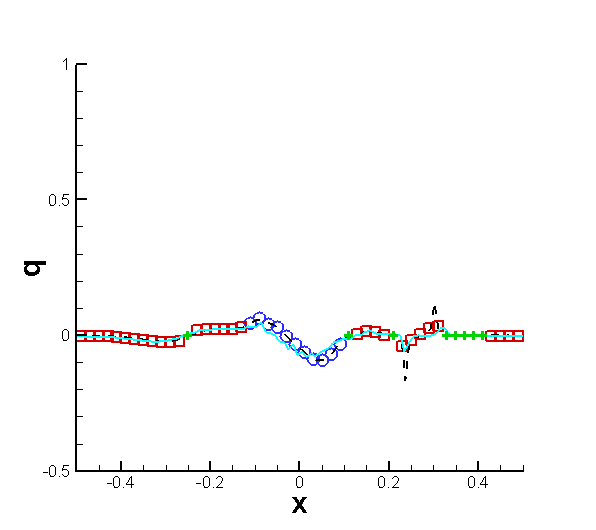}
\caption{Mixed regime problem with $\eps(x)$ with $a_0=40$ on the domain $[-0.5,0.5]\times[-10,10]$. $N_x=50$ and $N_v=100$ with NDG3.
Symbols: `+' is Euler, square is NS, circle is kinetic. Solid line: the reference NS solution with $N_x=100$ and $N_v=100$. Dashed line: the reference kinetic solution with $N_x=100$ and $N_v=100$. From left to right: time $t=0.1, 0.3, 0.45$. From top to bottom, the density $\rho$, the mean velocity $u$, the temperature $T$ and the heat flux $q$. $\eps_0=10^{-3}$, Euler-NS-kinetic.}
\label{fig53}
\end{figure}


\section{Conclusion}
\label{seccon}
\setcounter{equation}{0}
\setcounter{figure}{0}
\setcounter{table}{0}

We propose a high order hierarchical DG solver for the multi-scale BGK equation. Such hierarchical solver is based on an asymptotic preserving DG IMEX scheme \cite{JLQX_BGK}, which is formally showed to become a DG scheme for the limiting Euler system and a local DG scheme for the Navier-Stokes system when the Knudsen number is small. Adaptive criteria \cite{Filbet_hybrid} are applied to automatically switch DG solvers among different regimes (Euler, Navier-Stokes and kinetics), well balancing computational effectiveness and efficiency. Extensive numerical experiments are performed to showcase the proposed scheme in its ability for capturing solution structures and in computational savings.

\bibliographystyle{siam}
\bibliography{refer}

\end{document}